\documentclass[a4paper,12pt]{amsart}
\usepackage{amsmath}
\usepackage{amssymb}
\usepackage{amsthm,amsxtra}
\usepackage{latexsym, mathrsfs, color, soul}
\usepackage{cite}
\usepackage{enumerate}
\usepackage{enumitem}
\usepackage[mathscr]{eucal}
\usepackage{multicol}
\usepackage{multirow}
\usepackage{adjustbox}
\usepackage{float}
\usepackage{graphicx}
\usepackage{mathtools}
\usepackage[hyphens]{url}
\usepackage[colorlinks=true]{hyperref}
\usepackage{cleveref}
\usepackage{tikz-cd}
\usepackage{tkz-graph}
\usepackage{tkz-berge}
\usetikzlibrary{positioning,arrows,shapes.geometric,trees,decorations.text}
\usepackage[T1]{fontenc}
\usepackage{lmodern}
\usepackage[top=2in, bottom=1.5in, left=1in, right=1in]{geometry}

\makeatletter
\renewcommand{\@makefnmark}{%
  \hbox{\kern0.15em\@textsuperscript{\normalfont\@thefnmark}}%
}
\makeatother

\newtheorem{Def}{Definition}[section]
\newtheorem{Th}[Def]{Theorem}
\newtheorem{Ex}[Def]{Example}
\newtheorem{Lemma}[Def]{Lemma}
\newtheorem{Prop}[Def]{Proposition}
\newtheorem{Cor}[Def]{Corollary}
\newtheorem{Rem}[Def]{Remark}
\newtheorem{Obs}[Def]{Observation}

{}
\newtheorem{Prob}[Def]{Problem}

\DeclareMathOperator{\Uf}{U_{fin}}

\DeclareMathOperator{\S1k}{S_1^\kappa}
\DeclareMathOperator{\Sk}{S_{< \kappa}^\kappa}
\DeclareMathOperator{\Uk}{U_{< \kappa}^\kappa}

\DeclareMathOperator{\bk}{\mathfrak{b}_\kappa}
\DeclareMathOperator{\dk}{\mathfrak{d}_\kappa}

\DeclareMathOperator{\dlk}{\mathfrak{d}_\kappa^{<\kappa}}

\DeclareMathOperator{\maxk}{\max\nolimits_{<\kappa}}
\DeclareMathOperator{\lessk}{\mathord{<}\kappa}
\DeclareMathOperator{\ink}{in(\kappa)}

\DeclareMathOperator{\cf}{cf}
\DeclareMathOperator{\cof}{cof}

\DeclareMathOperator{\dom}{dom}

\DeclareMathAlphabet{\mybb}{U}{BOONDOX-ds}{m}{n}

\begin{document}
\title[Star selection principles in an uncountable setting]{Star selection principles in an uncountable setting}
\author[ N. Alam, D. Chandra ]{ Nur Alam$^\dag$, Debraj Chandra$^*$ }
\newcommand{\acr}{\newline\indent}
\address{\llap{\dag\,}Department of Mathematics, Mirza Ghalib College, Gaya, Bihar-823001, India}
\email{nurrejwana@gmail.com}
\address{\llap{*\,}Department of Mathematics, University of Gour Banga, Malda-732103, West Bengal, India}
\email{debrajchandra1986@gmail.com}


\maketitle

\begin{abstract}
The study of selection principles in an uncountable setting was initiated by Korch and Weiss (2020) and subsequently developed in \cite{ucsp, AC}. In the countable setting, star selection principles, introduced by Ko\v{c}inac (1999) as a generalization of selection principles via the star operation, have been extensively studied, particularly in $\Psi$-spaces and Pixley--Roy hyperspaces, and continue to be an active area of research. As a natural extension of this line of research, we initiate the study of star selection principles in the uncountable setting by introducing the star and strongly star variants of certain uncountable selection properties, namely the $\kappa$-Menger, $\kappa$-Scheepers, $\kappa$-Hurewicz, and $\kappa$-Rothberger properties. As a first and fundamental step, we investigate the relationships among these properties, establish equivalences under suitable topological assumptions, and provide counterexamples distinguishing them. Following the developments in the countable setting, a substantial part of our investigation is devoted to generalized $\Psi$-spaces, which provide a natural uncountable analogue of the classical $\Psi$-spaces. 
Several of our results on generalized $\Psi$-spaces reveal contrasting behaviour of star selection properties between the countable and uncountable settings.
\end{abstract}

\smallskip

\noindent{\bf\keywordsname{}:} {Uncountable selection principles, uncountable star selection principles, (strongly) star-$\kappa$-Menger, (strongly) star-$\kappa$-Scheepers, (strongly) star-$\kappa$-Hurewicz, (strongly) star-$\kappa$-Rothberger, generalized $\Psi$-spaces, $\mathfrak{d}_\kappa$, $\mathfrak{d}_\kappa^{<\kappa}$, $\mathfrak{b}_\kappa$, $\mathrm{in}(\kappa)$.}

\section{Definitions and Terminology}
\label{sec-dt}
Throughout, unless otherwise stated, $\kappa$ denotes an uncountable cardinal, and by a space we mean a topological space. For undefined notions and terminology, see \cite{AC, ucsp, kw20, engel, jech03, hsz26, hs25}.

\subsection{Uncountable selection principles}
\label{sec-csp}
We recall the following generalized notions of covers and selection principles from \cite{ucsp} (see also \cite{kw20, AC}). Let $\mathcal{U}$ be a cover of a space $X$. The family $\mathcal{U}$ is called a $\kappa$-cover of $X$ if $X\notin \mathcal{U}$ and for every $F\in [X]^{<\kappa}$, there exists $U\in \mathcal{U}$ such that $F\subseteq U$. The family $\mathcal{U}$ is called a $\gamma_\kappa$-cover of $X$ if $|\mathcal{U}| = \kappa$ and for every $x\in X$, the set $\{U\in \mathcal{U}: x\notin U\}$ has cardinality less than $\kappa$. Let $\mathcal{O}_\kappa$, $\Omega_\kappa$, and $\Gamma_\kappa$ denote, respectively, the classes of all open covers of cardinality $\kappa$, all open $\kappa$-covers of cardinality $\kappa$, and all open $\gamma_\kappa$-covers. For a class $\mathcal{A}$ of covers of spaces, let $\mathcal{A}(X)$ denote the family of all covers of $X$ belonging to $\mathcal{A}$. For a regular cardinal $\kappa$, the inclusions $\Gamma_\kappa(X) \subseteq \Omega_\kappa(X) \subseteq \mathcal{O}_\kappa(X)$ hold.

Let $\mathcal{A}, \mathcal{B}\in \{\mathcal{O}_\kappa, \Omega_\kappa, \Gamma_\kappa\}$. The uncountable selection principles for a space $X$ are formulated as follows.

\begin{enumerate}[leftmargin=2cm,itemsep=0.3em]
\item[$\S1k(\mathcal{A},\mathcal{B})$:] For each $(\mathcal{U}_\alpha : \alpha< \kappa) \in \mathcal{A}(X)^\kappa$, there are $U_\alpha\in \mathcal{U}_\alpha$, $\alpha< \kappa$, such that $\{U_\alpha : \alpha< \kappa\}\in \mathcal{B}(X)$.

\item[$\Sk(\mathcal{A},\mathcal{B})$:] For each $(\mathcal{U}_\alpha : \alpha< \kappa) \in \mathcal{A}(X)^\kappa$, there are $\mathcal{V}_\alpha\in [\mathcal{U}_\alpha]^{<\kappa}$, $\alpha< \kappa$, such that $\bigcup_{\alpha< \kappa} \mathcal{V}_\alpha \in \mathcal{B}(X)$.

\item[$\Uk(\mathcal{A},\mathcal{B})$:] For each $(\mathcal{U}_\alpha : \alpha< \kappa) \in \mathcal{A}(X)^\kappa$ with no $\mathcal{U}_\alpha$ has a subcover of size less than $\kappa$, there are $\mathcal{V}_\alpha\in [\mathcal{U}_\alpha]^{<\kappa}$, $\alpha< \kappa$, such that $\{\bigcup \mathcal{V}_\alpha : \alpha< \kappa\}\in \mathcal{B}(X)$.
\end{enumerate}

A space $X$ is said to be $\kappa$-Menger (respectively, $\kappa$-Scheepers, $\kappa$-Hurewicz, and $\kappa$-Rothberger) if it satisfies $\Sk(\mathcal{O}_\kappa, \mathcal{O}_\kappa)$ (respectively, $\Uk(\mathcal{O}_\kappa, \Omega_\kappa)$, $\Uk(\mathcal{O}_\kappa, \Gamma_\kappa)$, and $\S1k(\mathcal{O}_\kappa, \mathcal{O}_\kappa)$); see \cite{ucsp, kw20, AC}. A space $X$ is said to be $\kappa$-compact (respectively, $\kappa$-Lindel\"{o}f) if every open cover of $X$ has a subcover of cardinality less than $\kappa$ (respectively, at most $\kappa$). A family $\mathcal{U}$ of subsets of a space $X$ is said to be point-$\lessk$ if the size of $\{U\in \mathcal{U} : x\in U\}$ is less than $\kappa$ for every $x\in X$. Recall from \cite{fg16} that a space $X$ is said to be meta-$\lessk$-compact if every open cover of $X$ admits a point-$\lessk$ open refinement. Also recall that a space $X$ is said to be $\kappa$-paracompact if every open cover of $X$ of cardinality at most $\kappa$ admits a locally finite open refinement.

\subsection{Generalized Cantor and Baire spaces and cardinal invariants}
\label{sec-cb}
For a regular cardinal $\kappa$, consider a family $(X_\alpha : \alpha< \kappa)$ of discrete spaces. The $\kappa$-product topology on $\prod_{\alpha< \kappa} X_\alpha$ is the topology generated by the sets $\prod_{\alpha< \kappa} U_\alpha$, where $U_\alpha \subseteq X_\alpha$ for every $\alpha< \kappa$ and the size of $\{\alpha< \kappa : U_\alpha \neq X_\alpha\}$ is less than $\kappa$; see \cite{ucsp}. The spaces $2^\kappa$ and $\kappa^\kappa$, endowed with the $\kappa$-product topology, are called the $\kappa$-Cantor space and the $\kappa$-Baire space, respectively. By identifying each subset of $\kappa$ with its characteristic function, we regard $P(\kappa)$ as a topological space homeomorphic to $2^\kappa$. On the other hand, each $a\in [\kappa]^\kappa$ can be identified with its increasing enumeration in $\kappa^\kappa$. The two topologies on $[\kappa]^\kappa$ arising from these identifications coincide.

Identifying the elements of $[\kappa]^\kappa$ with their increasing enumerations, we write $a\leq_\kappa^* b$ for $a,b\in [\kappa]^\kappa$ if the size of $\{\alpha< \kappa : a(\alpha)> b(\alpha)\}$ is less than $\kappa$. A family $A\subseteq [\kappa]^\kappa$ is said to be bounded if there exists $b\in [\kappa]^\kappa$ such that $a\leq_\kappa^* b$ for every $a\in A$. The generalized bounding number $\bk$ is the least cardinality of an unbounded family in $[\kappa]^\kappa$. A family $A\subseteq [\kappa]^\kappa$ is said to be dominating if for every $x\in [\kappa]^\kappa$, there exists $a\in A$ such that $x\leq_\kappa^* a$. The generalized dominating number $\dk$ is the least cardinality of a dominating family in $[\kappa]^\kappa$. The notions of boundedness and domination can be defined analogously in $\kappa^\kappa$, yielding the same cardinal invariants $\bk$ and $\dk$.

For a regular cardinal $\kappa$ and a family $A\subseteq \kappa^\kappa$, put \[ \maxk(A) = \{\alpha \mapsto \sup \{f(\alpha) : f\in F\} : F\in [A]^{<\kappa}\},\] where $\sup \emptyset = 0$. Following \cite[Definition~1.11]{mrs05}, a family $A\subseteq \kappa^\kappa$ is said to be $\lessk$-dominating if for every $g\in \kappa^\kappa$, there exists $F\in [A]^{<\kappa}$ such that $g(\alpha)< \sup \{f(\alpha) : f\in F\}$ for every $\alpha< \kappa$. The $\lessk$-dominating number is defined by \[\dlk = \min \{|A| : A\subseteq \kappa^\kappa \text{ is $\lessk$-dominating}\}.\] We recall the following two lemmas from \cite{AC}.

\begin{Lemma}[{\cite{AC}}]
\label{lemma5}
Assume that $\kappa$ is regular. Then \[\dlk = \min \{|A| : A\subseteq \kappa^\kappa \text{ and } \maxk(A) \text{ is dominating}\}.\]
\end{Lemma}

\begin{Lemma}[{\cite{AC}}]
\label{lemma9}
Assume that $\kappa$ is regular. Then $\bk\leq \dlk\leq \dk$.
\end{Lemma}

Following Khomskii et al.~\cite[Definition~2.15]{klls16}, for a regular cardinal $\kappa$, two functions $f,g\in \kappa^\kappa$ are said to be eventually different if there exists $\beta< \kappa$ such that $f(\alpha) \neq g(\alpha)$ for every $\alpha \geq \beta$; equivalently, the size of $\{\alpha< \kappa : f(\alpha) = g(\alpha)\}$ is less than $\kappa$. A family $E\subseteq \kappa^\kappa$ is said to be eventually different if for every $g\in \kappa^\kappa$, there exists $f\in E$ that is eventually different from $g$. The inequality number is defined by \[\ink = \min \{|E| : E\subseteq \kappa^\kappa \text{ is an eventually different family}\}.\] It is immediate that $\ink\leq \dk$.

A cardinal $\kappa$ is weakly compact if and only if $2^\kappa$ is $\kappa$-compact \cite[Theorem~8.23]{cn74}. Every weakly compact cardinal $\kappa$ is strongly inaccessible; that is, $\kappa$ is regular and $2^\alpha< \kappa$ for every $\alpha< \kappa$. In particular, $\kappa^{<\kappa} = \kappa$. For an ordinal $\alpha$, let $\cf(\alpha)$ denote its cofinality. For a partially ordered set $(P,\leq)$, let $\cof(P,\leq)$, or simply $\cof(P)$ when the underlying order is clear from the context, denote its cofinality.

\section{Background and Motivation}
\label{sec-bm}
In \cite{coc1} (see also \cite{coc2}), Scheepers initiated the systematic study of selection principles in general and set-theoretic topology in the classical countable setting. The investigation of their uncountable analogues was initiated by Korch and Weiss in \cite{kw20} and subsequently developed in \cite{ucsp, AC}.

In \cite{kocinac99}, Ko\v{c}inac introduced and studied star selection principles, a generalization of the classical selection principles based on the star operation. For a subset $A$ of a space $X$ and a collection $\mathcal{P}$ of subsets of $X$, the star of $A$ with respect to $\mathcal{P}$ \cite{engel}, denoted by $St(A,\mathcal{P})$, is defined by \[St(A, \mathcal{P}) = \bigcup \{B\in \mathcal{P} : A\cap B\neq \emptyset\}.\] For $x\in X$, we write $St(x, \mathcal{P})$ instead of $St(\{x\}, \mathcal{P})$.

Since their introduction, star selection principles have been extensively investigated in the countable setting and continue to attract considerable attention in general and set-theoretic topology; see, for example, \cite{ac22, ac27, RRP, ca23, ca24, ca24-2, tsaban15, BCS, RBS, sx19, sx22} and the references therein. Motivated by the extensive study of star selection principles in the countable setting and the developments in uncountable selection principles \cite{kw20, ucsp, AC}, in this paper, we initiate the study of star selection principles in the uncountable setting by introducing the star variants of the $\kappa$-Menger, $\kappa$-Scheepers, $\kappa$-Hurewicz, and $\kappa$-Rothberger properties, namely, the star-$\kappa$-$P$ and strongly star-$\kappa$-$P$ properties, where $P\in$ \{Menger, Scheepers, Hurewicz, Rothberger\}.

We now introduce the star versions of the $\kappa$-Menger, $\kappa$-Scheepers, $\kappa$-Hurewicz, and $\kappa$-Rothberger properties as follows.

\begin{Def}
\label{def1}
A space $X$ is said to be
\begin{enumerate}[wide=0pt, label={\upshape(\arabic*)}, leftmargin=*, itemsep=0.3em]
\item star-$\kappa$-Menger (respectively, strongly star-$\kappa$-Menger) if for each $(\mathcal{U}_\alpha : \alpha< \kappa)\in \mathcal{O}_\kappa(X)^\kappa$, there are $\mathcal{V}_\alpha\in [\mathcal{U}_\alpha]^{<\kappa}$, $\alpha< \kappa$ (respectively, $F_\alpha\in [X]^{<\kappa}$, $\alpha< \kappa$), such that $X = \bigcup_{\alpha< \kappa} St(\bigcup \mathcal{V}_\alpha, \mathcal{U}_\alpha)$ (respectively, $X = \bigcup_{\alpha< \kappa} St(F_\alpha, \mathcal{U}_\alpha)$).

\item star-$\kappa$-Scheepers (respectively, strongly star-$\kappa$-Scheepers) if for each $(\mathcal{U}_\alpha : \alpha< \kappa)\in \mathcal{O}_\kappa(X)^\kappa$, there are $\mathcal{V}_\alpha\in [\mathcal{U}_\alpha]^{<\kappa}$, $\alpha< \kappa$ (respectively, $F_\alpha\in [X]^{<\kappa}$, $\alpha< \kappa$), such that for each $F\in [X]^{<\kappa}$, there exists an $\alpha< \kappa$ such that $F\subseteq St(\bigcup \mathcal{V}_\alpha, \mathcal{U}_\alpha)$ (respectively, $F\subseteq St(F_\alpha, \mathcal{U}_\alpha)$).

\item star-$\kappa$-Hurewicz (respectively, strongly star-$\kappa$-Hurewicz) if for each $(\mathcal{U}_\alpha : \alpha< \kappa)\in \mathcal{O}_\kappa(X)^\kappa$, there are $\mathcal{V}_\alpha\in [\mathcal{U}_\alpha]^{<\kappa}$, $\alpha< \kappa$ (respectively, $F_\alpha\in [X]^{<\kappa}$, $\alpha< \kappa$), such that for each $x\in X$, the set $\{\alpha< \kappa : x\notin St(\bigcup \mathcal{V}_\alpha, \mathcal{U}_\alpha)\}$ (respectively, $\{\alpha< \kappa : x\notin St(F_\alpha, \mathcal{U}_\alpha)\}$) has size less than $\kappa$.

\item star-$\kappa$-Rothberger (respectively, strongly star-$\kappa$-Rothberger) if for each $(\mathcal{U}_\alpha : \alpha< \kappa)\in \mathcal{O}_\kappa(X)^\kappa$, there are $U_\alpha\in \mathcal{U}_\alpha$, $\alpha< \kappa$ (respectively, $x_\alpha\in X$, $\alpha< \kappa$), such that $X = \bigcup_{\alpha< \kappa} St(U_\alpha, \mathcal{U}_\alpha)$ (respectively, $X = \bigcup_{\alpha< \kappa} St(x_\alpha, \mathcal{U}_\alpha)$).
\end{enumerate}
\end{Def}

Similarly, $X$ is said to be star-$\kappa$-compact (respectively, strongly star-$\kappa$-compact) if for every open cover $\mathcal{U}$ of $X$, there exists $\mathcal{V}\in [\mathcal{U}]^{<\kappa}$ (respectively, $F\in [X]^{<\kappa}$) such that $X = St(\bigcup \mathcal{V}, \mathcal{U})$ (respectively, $X = St(F,\mathcal{U})$).

In Section~\ref{sic}, we investigate the relationships among the $\kappa$-Menger, $\kappa$-Scheepers, $\kappa$-Hurewicz, and $\kappa$-Rothberger properties and their star variants. The implications among these properties are summarized in Figure~\ref{fig1}. We construct examples showing that none of the reverse implications in the figure holds in general, except for those along the third row, whose failure is established in \cite{AC}. These non-implications are summarized in Figure~\ref{fig2}. We further establish several equivalences under suitable topological assumptions. In particular, we show that the $\kappa$-Menger, $\kappa$-Scheepers, $\kappa$-Hurewicz, and $\kappa$-Rothberger properties are equivalent to their corresponding star variants in normal $\kappa$-paracompact spaces.

In the countable setting, $\Psi$-spaces \cite{mrowka54} have played a central role in the development of star selection principles, serving both as objects of extensive combinatorial investigation and as a rich source of counterexamples (see, for example, \cite{tsaban15, bm09, RBS, ca23, ca24}). In view of their significance, we consider generalized $\Psi$-spaces, that is, $\Psi_\kappa$-spaces \cite{xss20} as a natural starting point for our investigation of star selection principles in the uncountable setting. In Section~\ref{sca}, we establish combinatorial characterizations of the star-$\kappa$-Menger, star-$\kappa$-Scheepers, star-$\kappa$-Hurewicz, and star-$\kappa$-Rothberger properties, extending the corresponding results from the countable setting (see \cite[Theorems~2.1, 2.2, and 4.3]{tsaban15} and \cite[Theorem~3.12]{ca23}). We further establish cardinality characterizations of the $\kappa$-Menger, $\kappa$-Scheepers, and $\kappa$-Hurewicz properties in generalized $\Psi$-spaces. More precisely, for an almost disjoint family $\mathcal{A}\subseteq[\kappa]^\kappa$, we show that $\Psi_\kappa(\mathcal{A})$ is $\kappa$-Menger (respectively, $\kappa$-Scheepers, $\kappa$-Hurewicz) if and only if it is strongly star-$\kappa$-Menger (respectively, strongly star-$\kappa$-Scheepers, strongly star-$\kappa$-Hurewicz), which is equivalent to $|\mathcal{A}| < \dk$ (respectively, $|\mathcal{A}| < \dlk$, $|\mathcal{A}| < \bk$). These equivalences contrast with the corresponding situation in the classical countable setting. Since the $\kappa$-Menger, $\kappa$-Scheepers, and $\kappa$-Hurewicz properties are equivalent to their corresponding strongly star variants in meta-$\lessk$-compact as well as $\kappa$-paracompact spaces, we further show that the equivalences obtained for generalized $\Psi$-spaces are not restricted to the classes of meta-$\lessk$-compact or $\kappa$-paracompact spaces. In particular, assuming $\kappa^+ < \bk$, we construct a generalized $\Psi$-space that is strongly star-$\kappa$-Hurewicz but not $\kappa$-paracompact. Unlike the cases of the $\kappa$-Menger, $\kappa$-Scheepers, and $\kappa$-Hurewicz properties, whether the $\kappa$-Rothberger property is equivalent to its strongly star variant in generalized $\Psi$-spaces remains an open question (see Remark~\ref{rem1}). Nevertheless, in contrast to the classical countable setting, we also obtain a cardinality characterization of the $\kappa$-Rothberger property. More precisely, assuming $\kappa^{<\kappa} = \kappa$, we show that $\Psi_\kappa(\mathcal{A})$ is $\kappa$-Rothberger if and only if $|\mathcal{A}| < \ink$. Under the same assumption, we construct an almost disjoint family $\mathcal{A}$ of cardinality $\ink$ such that $\Psi_\kappa(\mathcal{A})$ fails to be star-$\kappa$-Rothberger, demonstrating that this cardinality bound is sharp for the strongly star-$\kappa$-Rothberger property.

Finally, we investigate the failure of star-$\kappa$ selection properties in generalized $\Psi$-spaces and obtain several contrasting results. In particular, we show that the natural uncountable analogues of certain results concerning star selection properties of $\Psi$-spaces in the countable setting do not hold in general (see \cite[Theorems~1.9 and 2.4]{tsaban15}). Assuming $\kappa^{<\kappa} = \kappa$, we construct an almost disjoint family $\mathcal{A}$ of cardinality $\lambda \geq \dk$ such that $\cof([\lambda]^\kappa) = \lambda$, while $\Psi_\kappa(\mathcal{A})$ is star-$\kappa$-Hurewicz, and consequently star-$\kappa$-Scheepers and star-$\kappa$-Menger. Thus the cardinality bounds $|\mathcal{A}| \geq \dk$, $|\mathcal{A}| \geq \dlk$, and $|\mathcal{A}| \geq \bk$, even together with the condition $\cof([|\mathcal{A}|]^\kappa) = |\mathcal{A}|$, do not necessarily imply the failure of the corresponding star-$\kappa$ selection properties. On the other hand, we establish sufficient structural conditions for generalized $\Psi$-spaces to fail these properties. In particular, we show that for every cardinal $\lambda$ with $\dk \leq \lambda \leq 2^\kappa$ (respectively, $\dlk \leq \lambda \leq 2^\kappa$, $\bk \leq \lambda \leq 2^\kappa$), there exists an almost disjoint family $\mathcal{A}$ of cardinality $\lambda$ such that $\Psi_\kappa(\mathcal{A})$ fails to be star-$\kappa$-Menger (respectively, star-$\kappa$-Scheepers, star-$\kappa$-Hurewicz). These results demonstrate that, beyond the corresponding cardinal thresholds, the star selection properties of generalized $\Psi$-spaces cannot be determined solely by the cardinality of the underlying almost disjoint family.

Similarly, for the star-$\kappa$-Rothberger property, we establish a sufficient combinatorial condition for its failure and show that for every cardinal $\lambda$ with $\ink \leq \lambda \leq 2^\kappa$, there exists an almost disjoint family $\mathcal{A}$ of cardinality $\lambda$ such that $\Psi_\kappa(\mathcal{A})$ is not star-$\kappa$-Rothberger. However, unlike the cases of the star-$\kappa$-Menger, star-$\kappa$-Scheepers, and star-$\kappa$-Hurewicz properties, we do not know whether there exists an almost disjoint family $\mathcal{A}$ with $|\mathcal{A}| \geq \ink$ and $\cof([|\mathcal{A}|]^\kappa) = |\mathcal{A}|$ such that $\Psi_\kappa(\mathcal{A})$ is star-$\kappa$-Rothberger. We leave this question open (see Problem~\ref{prob1}).

\section{Interconnections}
\label{sic}
\subsection{Implications and equivalences}
\label{sie}
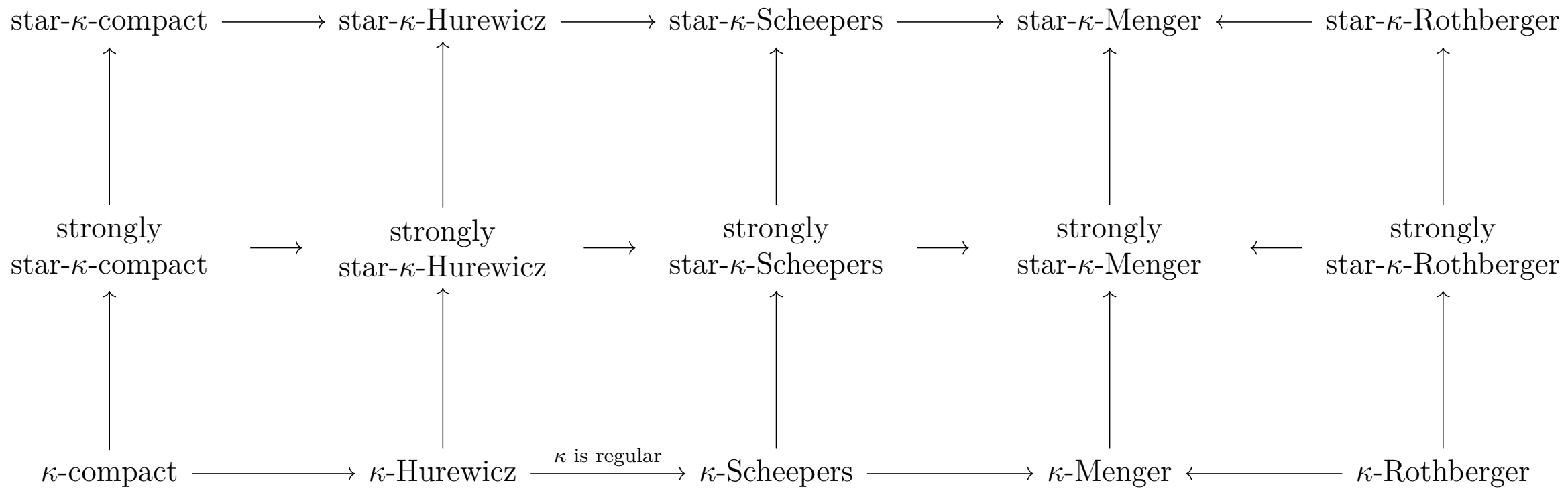
\begin{figure}[h]
\begin{adjustbox}{max width=\textwidth,max
height=\textheight,keepaspectratio,center}
\begin{tikzcd}[column sep=4ex, row sep=12ex, arrows={crossing over}]
\text{star-$\kappa$-compact} \arrow[r] & \text{star-$\kappa$-Hurewicz} \arrow[r] & \text{star-$\kappa$-Scheepers} \arrow[r] & \text{star-$\kappa$-Menger} & \text{star-$\kappa$-Rothberger} \arrow[l] \\
\text{\parbox{3.5cm}{\centering strongly\\ star-$\kappa$-compact}} \arrow[r] \arrow[u] & \text{\parbox{3.5cm}{\centering strongly\\ star-$\kappa$-Hurewicz}} \arrow[r] \arrow[u] & \text{\parbox{3.5cm}{\centering strongly\\ star-$\kappa$-Scheepers}} \arrow[r] \arrow[u] & \text{\parbox{3.5cm}{\centering strongly\\ star-$\kappa$-Menger}} \arrow[u] & \text{\parbox{3.5cm}{\centering strongly\\ star-$\kappa$-Rothberger}} \arrow[l] \arrow[u] \\
\text{$\kappa$-compact} \arrow[r]\arrow[u] & \text{$\kappa$-Hurewicz} \arrow[r,"\text{$\kappa$ is regular}"]\arrow[u] & \text{$\kappa$-Scheepers} \arrow[r]\arrow[u] & \text{$\kappa$-Menger} \arrow[u] & \text{$\kappa$-Rothberger} \arrow[l] \arrow[u]
\end{tikzcd}
\end{adjustbox}
\caption{Implications among $\kappa$-covering properties and their star variants.}
\label{fig1}
\end{figure}

The implications among the $\kappa$-compact, $\kappa$-Hurewicz, $\kappa$-Scheepers, $\kappa$-Menger, and $\kappa$-Rothberger properties and their star variants are summarized in Figure~\ref{fig1}.

Note that all the star variants depicted in Figure~\ref{fig1} are preserved under continuous images and inherited by clopen subsets.

\begin{Lemma}
\label{lemma26}
Let $\mathcal{U} \in \mathcal{O}_\kappa(X)$ have no subcover of cardinality less than $\kappa$. If $\mathcal{U}$ admits a point-$\lessk$ open refinement, then it admits a point-$\lessk$ open refinement of cardinality $\kappa$.
\end{Lemma}

\begin{Prop}
\label{prop1}
Let $P\in \{\text{Menger}, \text{Scheepers}, \text{Hurewicz}\}$. Assume that $\kappa$ is regular. Every meta-$\lessk$-compact (or $\kappa$-paracompact) strongly star-$\kappa$-$P$ space is $\kappa$-$P$.
\end{Prop}

\begin{proof}
We prove only the case when $P = \text{Scheepers}$ and the underlying space is meta-$\lessk$-compact, as the other cases follow by analogous arguments. Let $X$ be a meta-$\lessk$-compact strongly star-$\kappa$-Scheepers space. Pick $(\mathcal{U}_\alpha : \alpha< \kappa) \in \mathcal{O}_\kappa(X)^\kappa$ such that no $\mathcal{U}_\alpha$ has a subcover of cardinality less than $\kappa$. For each $\alpha< \kappa$, choose a point-$\lessk$ open refinement $\mathcal{V}_\alpha$ of $\mathcal{U}_\alpha$. By Lemma~\ref{lemma26}, we can assume that $|\mathcal{V}_\alpha| = \kappa$ for every $\alpha< \kappa$. Since $X$ is strongly star-$\kappa$-Scheepers, there exist $F_\alpha \in [X]^{<\kappa}$, $\alpha< \kappa$, such that for every $A\in [X]^{<\kappa}$, there exists $\alpha< \kappa$ satisfying $A\subseteq St(F_\alpha, \mathcal{V}_\alpha)$. For each $\alpha< \kappa$, put $\mathcal{W}_\alpha = \{V\in \mathcal{V}_\alpha : V\cap F_\alpha \neq \emptyset\}$. Since $\mathcal{V}_\alpha$ is point-$\lessk$, $|F_\alpha|< \kappa$, and $\kappa$ is regular, we have $|\mathcal{W}_\alpha| < \kappa$. For every $V\in \mathcal{W}_\alpha$, choose $U_V\in \mathcal{U}_\alpha$ such that $V\subseteq U_V$, and put $\mathcal{H}_\alpha = \{U_V : V\in \mathcal{W}_\alpha\}$. Then
$\mathcal{H}_\alpha \in [\mathcal{U}_\alpha]^{<\kappa}$. Moreover, for every $A\in [X]^{<\kappa}$, there exists $\alpha< \kappa$ such that $A\subseteq St(F_\alpha, \mathcal{V}_\alpha) = \bigcup \mathcal{W}_\alpha \subseteq \bigcup \mathcal{H}_\alpha$.

We now enlarge the selections so that their unions are pairwise distinct. For each $\alpha< \kappa$, let $\mathcal{D}_\alpha = \{\bigcup (\mathcal{H}_\alpha \cup \{U\}) : U\in \mathcal{U}_\alpha\}$. We claim that $|\mathcal{D}_\alpha| = \kappa$. Otherwise, choosing one representative $U_D\in \mathcal{U}_\alpha$ for each $D\in \mathcal{D}_\alpha$, we obtain $X = \bigcup \mathcal{H}_\alpha \cup \bigcup_{D\in \mathcal{D}_\alpha} U_D$. This yields a subcover of $\mathcal{U}_\alpha$ of cardinality less than $\kappa$, a contradiction. Thus $|\mathcal{D}_\alpha| = \kappa$ for every $\alpha< \kappa$. By transfinite recursion, choose $U_\alpha\in\mathcal{U}_\alpha$, $\alpha< \kappa$, such that the sets $\bigcup (\mathcal{H}_\alpha \cup \{U_\alpha\})$, $\alpha< \kappa$, are pairwise distinct.

For each $\alpha< \kappa$, define $\mathcal{G}_\alpha = \mathcal{H}_\alpha \cup \{U_\alpha\}$. Clearly, $\mathcal{G}_\alpha \in [\mathcal{U}_\alpha]^{<\kappa}$. Since $\mathcal{U}_\alpha$ has no subcover of cardinality less than $\kappa$, we have $\bigcup \mathcal{G}_\alpha \neq X$ for every $\alpha< \kappa$. Furthermore, the sets $\bigcup \mathcal{G}_\alpha$ are pairwise distinct, and for every $A\in [X]^{<\kappa}$, there exists $\alpha< \kappa$ such that $A\subseteq \bigcup \mathcal{H}_\alpha \subseteq \bigcup \mathcal{G}_\alpha$. Therefore, $\{\bigcup \mathcal{G}_\alpha : \alpha< \kappa\} \in \Omega_\kappa(X)$. Hence $X$ is $\kappa$-Scheepers.
\end{proof}

\begin{Cor}
\label{cor8}
Let $P\in \{\text{Menger}, \text{Scheepers}, \text{Hurewicz}\}$. Assume that $\kappa$ is regular. Every meta-Lindel\"{o}f strongly star-$\kappa$-$P$ space is $\kappa$-$P$. In particular, every metacompact strongly star-$\kappa$-$P$ space is $\kappa$-$P$.
\end{Cor}

Recall from \cite{engel} that an open cover $\mathcal{V}$ of a space $X$ is said to be a point-star refinement of an open cover $\mathcal{U}$ of $X$ if for every $x\in X$, there exists $U\in \mathcal{U}$ such that $St(x,\mathcal{V}) \subseteq U$.

\begin{Lemma}
\label{lemma27}
Let $X$ be a normal $\kappa$-paracompact space. Then every open cover of $X$ of cardinality at most $\kappa$ admits an open star refinement of cardinality at most $\kappa$. Moreover, if the original cover has no subcover of cardinality less than $\kappa$, then the star refinement has cardinality exactly $\kappa$.
\end{Lemma}

\begin{proof}
Let $\mathcal{U}$ be an open cover of $X$ with $|\mathcal{U}| \leq \kappa$. By $\kappa$-paracompactness, there exists a locally finite open refinement $\mathcal{W}$ of $\mathcal{U}$. By grouping together the members of $\mathcal{W}$ according to their chosen containing members of $\mathcal{U}$, we obtain a locally finite open refinement of $\mathcal{U}$ of cardinality at most $\kappa$. Replacing $\mathcal{W}$ by this refinement, we may assume that $|\mathcal{W}| \leq \kappa$. Write $\mathcal{W} = \{W_i : i\in I\}$, $|I| \leq \kappa$. Since $X$ is normal and $\mathcal{W}$ is locally finite, there exist closed sets $F_i \subseteq W_i$, $i\in I$, such that $X = \bigcup_{i\in I} F_i$. For every nonempty finite $J\subseteq I$, put $V_J = \bigcap_{i\in J} W_i \setminus \bigcup_{i\in I\setminus J} F_i$. Since $\{F_i : i\in I\}$ is locally finite, $V_J$ is open. Let \[\mathcal{V} = \{V_J : J\in [I]^{<\omega} \setminus
\{\emptyset\}, V_J\neq \emptyset\}.\] It is easy to see that $\mathcal{V}$ is an open cover of $X$ and $|\mathcal{V}| \leq \kappa$. For every $x\in X$, choose $i\in I$ such that $x\in F_i$. If $x\in V_J$, then $i\in J$. Consequently, $St(x,\mathcal{V}) \subseteq W_i$. Thus $\mathcal{V}$ is a point-star refinement of $\mathcal{W}$.

Applying the same construction to $\mathcal{V}$, we obtain a point-star refinement $\mathcal{Z}$ of $\mathcal{V}$ with $|\mathcal{Z}| \leq \kappa$. Since a point-star refinement of a point-star refinement is a star refinement (see \cite[Lemma~5.1.15]{engel}), $\mathcal{Z}$ is a star refinement of $\mathcal{W}$, and hence of $\mathcal{U}$.

Finally, suppose that $\mathcal{U}$ has no subcover of cardinality less than $\kappa$. If $|\mathcal{Z}| < \kappa$, then choosing for each $Z\in \mathcal{Z}$, a member of $\mathcal{U}$ containing $Z$ would yield a subcover of $\mathcal{U}$ of cardinality less than $\kappa$, a contradiction. Therefore, $|\mathcal{Z}| = \kappa$.
\end{proof}

\begin{Lemma}
\label{lemma28}
A space $X$ is $\kappa$-Menger if and only if for every sequence $(\mathcal{U}_\alpha : \alpha< \kappa) \in \mathcal{O}_\kappa(X)^\kappa$ such that no $\mathcal{U}_\alpha$ has a subcover of cardinality less than $\kappa$, there exist $\mathcal{V}_\alpha \in [\mathcal{U}_\alpha]^{<\kappa}$, $\alpha< \kappa$, such that $\bigcup_{\alpha< \kappa} \mathcal{V}_\alpha$ covers $X$.
\end{Lemma}

\begin{Prop}
\label{prop4}
Let $P\in \{\text{Menger}, \text{Scheepers}, \text{Hurewicz}, \text{Rothberger}\}$. Assume that $\kappa$ is regular. Every normal $\kappa$-paracompact star-$\kappa$-$P$ space is $\kappa$-$P$.
\end{Prop}

\begin{proof}
It suffices to prove the case $P = \text{Rothberger}$, since the remaining cases can be established similarly. Let $X$ be a normal $\kappa$-paracompact star-$\kappa$-Rothberger space. Pick $(\mathcal{U}_\alpha : \alpha< \kappa)\in \mathcal{O}_\kappa(X)^\kappa$. First suppose that $w(X)< \kappa$. Let $\mathcal{B}$ be a base of $X$ with $|\mathcal{B}| < \kappa$. For each $B\in \mathcal{B}$, put $A_B = \{\alpha< \kappa : (\exists U\in \mathcal{U}_\alpha)\, B\subseteq U\}$. By the regularity of $\kappa$, the family $\mathcal{B}^\prime = \{B\in \mathcal{B} : |A_B| = \kappa\}$ covers $X$. Since $|\mathcal{B}^\prime| < \kappa$, we may choose pairwise distinct $\alpha_B\in A_B$, $B\in \mathcal{B}^\prime$. For each $B\in \mathcal{B}^\prime$, choose $U_{\alpha_B} \in \mathcal{U}_{\alpha_B}$ such that $B\subseteq U_{\alpha_B}$. These choices cover $X$. The choices at the remaining indices can be made pairwise distinct and distinct from the already chosen members, since each original cover has cardinality $\kappa$. Therefore, $X$ is $\kappa$-Rothberger.

Now assume that $w(X)\geq \kappa$. By Lemma~\ref{lemma27}, for every $\alpha< \kappa$, there exists an open star refinement $\mathcal{V}_\alpha$ of $\mathcal{U}_\alpha$ such that $|\mathcal{V}_\alpha| \leq \kappa$. We may assume that $|\mathcal{V}_\alpha| = \kappa$. Partition $\kappa$ into two disjoint sets $I$ and $J$, each of cardinality $\kappa$. Applying the star-$\kappa$-Rothberger property of $X$ to the sequence $(\mathcal{V}_\alpha : \alpha\in I)$, reindexed by $\kappa$, we obtain $V_\alpha\in \mathcal{V}_\alpha$, $\alpha\in I$, such that $X = \bigcup_{\alpha\in I} St(V_\alpha, \mathcal{V}_\alpha)$. For each $\alpha\in I$, choose $U_\alpha\in \mathcal{U}_\alpha$ satisfying $St(V_\alpha, \mathcal{V}_\alpha) \subseteq U_\alpha$. It follows that $X = \bigcup_{\alpha\in I} St(V_\alpha, \mathcal{V}_\alpha) \subseteq \bigcup_{\alpha\in I} U_\alpha$. For each $\alpha\in J$, recursively choose $U_\alpha \in \mathcal{U}_\alpha \setminus \{U_\beta : \beta\in J, \beta< \alpha\}$. Such a choice is possible since $|\mathcal{U}_\alpha| = \kappa$ and fewer than $\kappa$ members are excluded at each stage. Thus $\{U_\alpha : \alpha< \kappa\}\in \mathcal{O}_\kappa(X)$. Hence $X$ is $\kappa$-Rothberger.
\end{proof}

\begin{Cor}
\label{cor13}
Let $P\in \{\text{Menger}, \text{Scheepers}, \text{Hurewicz}, \text{Rothberger}\}$. Assume that $\kappa$ is regular. For a normal $\kappa$-paracompact (or paracompact Hausdorff) space $X$, the following assertions are equivalent.
\begin{enumerate}[wide=0pt, label={\upshape(\arabic*)}, leftmargin=*]
\item $X$ is star-$\kappa$-$P$.

\item $X$ is strongly star-$\kappa$-$P$.

\item $X$ is $\kappa$-$P$.
\end{enumerate}
\end{Cor}

\subsection{Non-implications}
\label{sni}
In this section, we show that the reverse implications in Figure~\ref{fig1} do not hold in general. For those along the third row, see \cite{AC}. To establish the remaining non-implications, we use the following lemmas. We provide proofs only for Lemmas~\labelcref{lemma4,lemma11}, as the proofs of the remaining lemmas are either analogous or elementary.

\begin{Lemma}
\label{lemma7}
If $\{X_\xi : \xi< \kappa\}$ is a family of star-$\kappa$-Menger (respectively, strongly star-$\kappa$-Menger) subspaces of a space $X$, then $Y = \bigcup_{\xi< \kappa} X_\xi$ is star-$\kappa$-Menger (respectively, strongly star-$\kappa$-Menger).
\end{Lemma}

\begin{Lemma}
\label{lemma6}
If $X = Y\cup Z$, where $Y$ is star-$\kappa$-Scheepers (respectively, strongly star-$\kappa$-Scheepers) and $Z$ is star-$\kappa$-compact (respectively, strongly star-$\kappa$-compact), then $X$ is star-$\kappa$-Scheepers (respectively, strongly star-$\kappa$-Scheepers).
\end{Lemma}

\begin{Lemma}
\label{lemma4}
If $\{X_\xi : \xi< \cf(\kappa)\}$ is a family of star-$\kappa$-Hurewicz (respectively, strongly star-$\kappa$-Hurewicz) subspaces of a space $X$, then $Y = \bigcup_{\xi< \cf(\kappa)} X_\xi$ is star-$\kappa$-Hurewicz (respectively, strongly star-$\kappa$-Hurewicz).
\end{Lemma}

\begin{proof}
We prove only the star-$\kappa$-Hurewicz case, as the proof of the other case is similar. Choose a strictly increasing sequence $(\gamma_\xi : \xi< \cf(\kappa))$ of ordinals below $\kappa$ which is cofinal in $\kappa$. Let $(\mathcal{U}_\alpha : \alpha< \kappa)\in \mathcal{O}_\kappa(Y)^\kappa$. For each $\xi< \cf(\kappa)$ and $\alpha< \kappa$, consider the open cover $\mathcal{U}_\alpha^{(\xi)} = \{U\cap X_\xi : U\in \mathcal{U}_\alpha, U\cap X_\xi \neq \emptyset\}$ of $X_\xi$. We first show that, for each $\xi< \cf(\kappa)$, there are $\mathcal{V}_\alpha^{(\xi)} \in [\mathcal{U}_\alpha]^{<\kappa}$, $\alpha< \kappa$, such that for each $x\in X_\xi$, the size of $\{\alpha< \kappa : x\notin St(\bigcup \mathcal{V}_\alpha^{(\xi)}, \mathcal{U}_\alpha)\}$ is less than $\kappa$.

Fix $\xi<\cf(\kappa)$, and put $J_\xi = \{\alpha< \kappa : |\mathcal{U}_\alpha^{(\xi)}| = \kappa\}$. If $J_\xi = \emptyset$, then for each
$\alpha< \kappa$, choose $\mathcal{V}_\alpha^{(\xi)} \in [\mathcal{U}_\alpha]^{<\kappa}$ such that $X_\xi \subseteq \bigcup \mathcal{V}_\alpha^{(\xi)}$. Clearly, for each $x\in X_\xi$, the size of $\{\alpha< \kappa : x\notin St(\bigcup \mathcal{V}_\alpha^{(\xi)}, \mathcal{U}_\alpha)\}$ is less than $\kappa$. Suppose now that $J_\xi \neq \emptyset$. Fix $\beta_\xi \in J_\xi$ and define \[\mathcal{W}_\alpha^{(\xi)} = \begin{cases}
\mathcal{U}_\alpha^{(\xi)}, & \text{if } \alpha \in J_\xi,\\
\mathcal{U}_{\beta_\xi}^{(\xi)}, & \text{if } \alpha \notin J_\xi.
\end{cases}\]
Then $(\mathcal{W}_\alpha^{(\xi)} : \alpha< \kappa)\in \mathcal{O}_\kappa(X_\xi)^\kappa$. Since $X_\xi$ is star-$\kappa$-Hurewicz,
there are $\mathcal{A}_\alpha^{(\xi)} \in [\mathcal{W}_\alpha^{(\xi)}]^{<\kappa}$, $\alpha< \kappa$, such that for each $x\in X_\xi$, the size of $\{\alpha< \kappa : x\notin St(\bigcup \mathcal{A}_\alpha^{(\xi)}, \mathcal{W}_\alpha^{(\xi)})\}$ is less than $\kappa$. For each $\alpha\in J_\xi$, choose $\mathcal{V}_\alpha^{(\xi)} \in [\mathcal{U}_\alpha]^{<\kappa}$ such that $\mathcal{A}_\alpha^{(\xi)} = \{U\cap X_\xi : U\in \mathcal{V}_\alpha^{(\xi)}\}$. For each $\alpha \notin J_\xi$, choose $\mathcal{V}_\alpha^{(\xi)} \in [\mathcal{U}_\alpha]^{<\kappa}$ such that $X_\xi \subseteq \bigcup \mathcal{V}_\alpha^{(\xi)}$. It follows that for each $x\in X_\xi$, the size of $\{\alpha< \kappa : x\notin St(\bigcup \mathcal{V}_\alpha^{(\xi)}, \mathcal{U}_\alpha)\}$ is less than $\kappa$.

For each $\alpha< \kappa$, put $I_\alpha = \{\xi< \cf(\kappa) : \gamma_\xi \leq \alpha\}$ and $\mathcal{V}_\alpha = \bigcup_{\xi\in I_\alpha} \mathcal{V}_\alpha^{(\xi)}$. Since $(\gamma_\xi : \xi< \cf(\kappa))$ is strictly increasing and cofinal in $\kappa$, we have $|I_\alpha|< \cf(\kappa)$. Consequently, $|\mathcal{V}_\alpha| < \kappa$. Thus $\mathcal{V}_\alpha \in [\mathcal{U}_\alpha]^{<\kappa}$ for every $\alpha< \kappa$. Let $x\in Y$. Choose $\xi< \cf(\kappa)$ such that $x\in X_\xi$. For every $\alpha \geq \gamma_\xi$, we have $\xi\in I_\alpha$, and hence $St(\bigcup \mathcal{V}_\alpha^{(\xi)}, \mathcal{U}_\alpha) \subseteq St(\bigcup \mathcal{V}_\alpha, \mathcal{U}_\alpha)$. Therefore, \[\{\alpha< \kappa : x\notin St(\bigcup \mathcal{V}_\alpha, \mathcal{U}_\alpha)\} \subseteq
\gamma_\xi \cup \{\alpha< \kappa : x\notin St(\bigcup \mathcal{V}_\alpha^{(\xi)}, \mathcal{U}_\alpha)\}.\] It follows that the size of $\{\alpha< \kappa : x\notin St(\bigcup \mathcal{V}_\alpha, \mathcal{U}_\alpha)\}$ is less than $\kappa$. Thus $Y$ is star-$\kappa$-Hurewicz.
\end{proof}

\begin{Lemma}
\label{lemma11}
If $X$ is a star-$\kappa$-Rothberger (respectively, strongly star-$\kappa$-Rothberger) space, then for every sequence $(\mathcal{U}_\alpha : \alpha< \kappa)$ of open covers of $X$ of size at most $\kappa$, there are $U_\alpha \in \mathcal{U}_\alpha$ (respectively, $x_\alpha \in X$), $\alpha< \kappa$, such that $X = \bigcup_{\alpha< \kappa} St(U_\alpha, \mathcal{U}_\alpha)$ (respectively, $X = \bigcup_{\alpha< \kappa} St(x_\alpha, \mathcal{U}_\alpha)$).
\end{Lemma}

\begin{proof}
We present the proof only for the star-$\kappa$-Rothberger case, since the other case is analogous. Let $(\mathcal{U}_\alpha : \alpha< \kappa)$ be a sequence of open covers of $X$ with $|\mathcal{U}_\alpha| \leq \kappa$ for every $\alpha< \kappa$. If every $\mathcal{U}_\alpha$ admits an open refinement of cardinality $\kappa$, then the desired conclusion follows from the star-$\kappa$-Rothberger property of $X$. Suppose that some $\mathcal{U}_{\alpha_0}$ has no open refinement of cardinality $\kappa$. Since $\mathcal{U}_{\alpha_0}$ is itself an open refinement, we have $|\mathcal{U}_{\alpha_0}| < \kappa$. Write $\mathcal{U}_{\alpha_0} = \{V_\xi : \xi< \mu\}$, where $\mu< \kappa$. For each $\xi< \mu$, there are fewer than $\kappa$ distinct open subsets of $X$ contained in $V_\xi$. Indeed, otherwise, adjoining $\kappa$ distinct open subsets of $V_\xi$ to $\mathcal{U}_{\alpha_0}$ would yield an open refinement of cardinality $\kappa$, a contradiction.

Partition $\kappa$ into pairwise disjoint sets $I_\xi$, $\xi< \mu$, each of
cardinality $\kappa$. Fix $\xi< \mu$. For each $\alpha\in I_\xi$, recursively choose $U_\alpha \in \mathcal{U}_\alpha$ as follows. If \[V_\xi \setminus \bigcup_{\substack{\beta\in I_\xi\\ \beta< \alpha}} U_\beta \neq \emptyset,\] choose a point $x_\alpha$ in this difference and take $U_\alpha \in \mathcal{U}_\alpha$ containing $x_\alpha$. Otherwise, choose $U_\alpha$ arbitrarily. We claim that $V_\xi \subseteq \bigcup_{\alpha\in I_\xi} U_\alpha$. Otherwise, for each $\alpha\in I_\xi$, \[V_\xi \setminus \bigcup_{\substack{\beta\in I_\xi\\ \beta< \alpha}} U_\beta \neq \emptyset.\] By construction, the sets $U_\alpha \cap V_\xi$, $\alpha\in I_\xi$, are $\kappa$ distinct open subsets of $X$ contained in $V_\xi$, a contradiction. Thus $V_\xi \subseteq \bigcup_{\alpha\in I_\xi} U_\alpha$ for every $\xi< \mu$. Consequently, $X = \bigcup_{\xi< \mu} V_\xi \subseteq \bigcup_{\alpha< \kappa} U_\alpha \subseteq \bigcup_{\alpha< \kappa} St(U_\alpha, \mathcal{U}_\alpha)$. Hence the result follows.
\end{proof}

As an immediate consequence of the preceding result, we obtain the following lemma.

\begin{Lemma}
\label{lemma10}
If $\{X_\xi : \xi< \kappa\}$ is a family of star-$\kappa$-Rothberger (respectively, strongly star-$\kappa$-Rothberger) subspaces of a space $X$, then $Y = \bigcup_{\xi< \kappa} X_\xi$ is star-$\kappa$-Rothberger (respectively, strongly star-$\kappa$-Rothberger).
\end{Lemma}

\begin{Lemma}
\label{lemma8}
Let $\lambda< \cf(\kappa)$. If $\{X_\xi : \xi< \lambda\}$ is a family of star-$\kappa$-compact (respectively, strongly star-$\kappa$-compact) subspaces of a space $X$, then $Y = \bigcup_{\xi< \lambda} X_\xi$ is star-$\kappa$-compact (respectively, strongly star-$\kappa$-compact).
\end{Lemma}

We begin by considering the reverse implications along the first row. The space considered in the following example is taken from \cite[Example~2.2]{song14}.

\begin{Ex}
\label{ex5}
For every cardinal $\kappa$, there exists a Tychonoff star-$\kappa$-Hurewicz space which is neither strongly star-$\kappa$-Menger nor star-$\kappa$-compact.
\end{Ex}

\begin{proof}
Let $\lambda = |\kappa^\kappa|$, and let $D = \{d_f : f\in \kappa^\kappa\}$ be a discrete space. Let $aD = D\cup \{d^*\}$ be the one-point compactification of $D$, and put \[Y(\kappa^\kappa) = (aD\times [0,\lambda^+)) \cup (D\times \{\lambda^+\}),\] considered as a subspace of $aD\times[0,\lambda^+]$.

We first show that $Y(\kappa^\kappa)$ is star-$\kappa$-compact. Let $\mathcal{U}$ be an open cover of $Y(\kappa^\kappa)$. For each $f\in \kappa^\kappa$, choose $U_f\in \mathcal{U}$ such that $(d_f, \lambda^+) \in U_f$. Since $U_f$ is open in $Y(\kappa^\kappa)$, there exists $\beta_f< \lambda^+$ such that $\{d_f\} \times (\beta_f, \lambda^+] \subseteq U_f$. Since $|\kappa^\kappa| = \lambda< \lambda^+ = \cf(\lambda^+)$, there exists $\beta< \lambda^+$ such that $\beta_f< \beta$ for every $f\in \kappa^\kappa$. The set $aD\times \{\beta\}$ is compact. Hence there exists a finite subfamily $\mathcal{V}_1 \subseteq \mathcal{U}$ such that $aD\times \{\beta\} \subseteq \bigcup \mathcal{V}_1$. For every $f\in \kappa^\kappa$, $(d_f, \beta)\in U_f\cap \bigcup \mathcal{V}_1$, and therefore $(d_f, \lambda^+) \in St \left(\bigcup \mathcal{V}_1, \mathcal{U}\right)$. Thus $D\times \{\lambda^+\} \subseteq St\left(\bigcup \mathcal{V}_1, \mathcal{U}\right)$. Since $aD$ is compact and $[0, \lambda^+)$ is countably compact, $aD\times [0,\lambda^+)$ is countably compact. Choose a finite subfamily $\mathcal{V}_2 \subseteq \mathcal{U}$ such that $aD\times [0,\lambda^+) \subseteq St\left(\bigcup \mathcal{V}_2, \mathcal{U}\right)$. Finally, put $\mathcal{V} = \mathcal{V}_1 \cup \mathcal{V}_2$. Then $\mathcal{V}$ is finite and $Y(\kappa^\kappa) = St\left(\bigcup \mathcal{V}, \mathcal{U}\right)$. Thus $Y(\kappa^\kappa)$ is star-$\kappa$-compact.

We now show that $Y(\kappa^\kappa)$ is not strongly star-$\kappa$-Menger. For
$\alpha, \xi< \kappa$, let $D^{(\alpha)}_\xi = \{d_f : f\in \kappa^\kappa, f(\alpha) = \xi\}$ and put $W^{(\alpha)}_\xi = (D^{(\alpha)}_\xi \times [0, \lambda^+]) \cap Y(\kappa^\kappa)$. Also let $V = aD\times [0, \lambda^+)$. Since $D^{(\alpha)}_\xi$ is open in $aD$, each $W^{(\alpha)}_\xi$ is open in $Y(\kappa^\kappa)$. Moreover, $V$ is open in $Y(\kappa^\kappa)$. For each $\alpha< \kappa$, define $\mathcal{U}_\alpha = \{V\} \cup \{W^{(\alpha)}_\xi : \xi< \kappa\}$. Then $\mathcal{U}_\alpha \in \mathcal{O}_\kappa(Y(\kappa^\kappa))$. Let $F_\alpha \in [Y(\kappa^\kappa)]^{<\kappa}$, $\alpha< \kappa$, be arbitrary. For each $\alpha< \kappa$, put \[ A_\alpha = \{\xi< \kappa : F_\alpha \cap W^{(\alpha)}_\xi \neq \emptyset\}.\] Since the sets $W^{(\alpha)}_\xi$, $\xi< \kappa$, are pairwise disjoint, $|A_\alpha| \leq |F_\alpha| < \kappa$. Choose $g(\alpha) \in \kappa \setminus A_\alpha$ for every $\alpha< \kappa$. Thus $g\in \kappa^\kappa$. Consider the point $x = (d_g, \lambda^+) \in Y(\kappa^\kappa)$. For each $\alpha< \kappa$, the only member of $\mathcal{U}_\alpha$ containing $x$ is $W^{(\alpha)}_{g(\alpha)}$. By the choice of $g(\alpha)$, $F_\alpha \cap W^{(\alpha)}_{g(\alpha)} = \emptyset$. Consequently, $x\notin St(F_\alpha, \mathcal{U}_\alpha)$ for every $\alpha< \kappa$. Hence $x\notin \bigcup_{\alpha< \kappa} St(F_\alpha, \mathcal{U}_\alpha)$. Therefore, $Y(\kappa^\kappa)$ is not strongly star-$\kappa$-Menger.

Let $Z = \{z_\xi : \xi< \kappa\}$ be a discrete space of cardinality $\kappa$, and put $X = Y(\kappa^\kappa) \oplus Z$. Clearly, $X$ is a Tychonoff space which is neither star-$\kappa$-compact nor strongly star-$\kappa$-Menger. Also by Lemma~\ref{lemma4}, $X$ is star-$\kappa$-Hurewicz.
\end{proof}

\begin{Ex}
\label{ex6}
Assume that $\kappa^{<\kappa} = \kappa$ and $\bk< \dk$. Then there exists a Tychonoff star-$\kappa$-Scheepers space which is neither star-$\kappa$-Hurewicz nor strongly star-$\kappa$-Menger.
\end{Ex}

\begin{proof}
Let $Y$ be the space  $Y(\kappa^\kappa)$ defined in Example~\ref{ex5}. Then $Y$ is a Tychonoff star-$\kappa$-compact space which is not strongly star-$\kappa$-Menger. Under the assumption $\kappa^{<\kappa} = \kappa$ and $\bk< \dk$, there exists a $\kappa$-Scheepers subspace $Z$ of $\kappa^\kappa$ which is not $\kappa$-Hurewicz (see \cite{AC}). Put $X = Y\oplus Z$. Clearly, $X$ is a Tychonoff space which is neither star-$\kappa$-Hurewicz nor strongly star-$\kappa$-Scheepers. By Lemma~\ref{lemma6}, $X$ is star-$\kappa$-Scheepers. \end{proof}

\begin{Ex}
\label{ex7}
Assume that $\kappa^{<\kappa} = \kappa$ and $\dlk< \dk$. Then there exists a Tychonoff star-$\kappa$-Menger space which is neither star-$\kappa$-Scheepers nor strongly star-$\kappa$-Menger.
\end{Ex}

\begin{proof}
Let $Y$ be the space  $Y(\kappa^\kappa)$ defined in Example~\ref{ex5}. Then $Y$ is a Tychonoff star-$\kappa$-compact space which is not strongly star-$\kappa$-Menger. Under the assumption $\kappa^{<\kappa} = \kappa$ and $\dlk< \dk$, there exists a $\kappa$-Menger subspace $Z$ of $\kappa^\kappa$ which is not $\kappa$-Scheepers (see \cite{AC}). Put $X = Y\oplus Z$. Clearly, $X$ is a Tychonoff space which is neither star-$\kappa$-Scheepers nor strongly star-$\kappa$-Menger. By Lemma~\ref{lemma7}, $X$ is star-$\kappa$-Menger.
\end{proof}

\begin{Ex}
\label{ex8}
Let $\kappa$ be a weakly compact cardinal. Then there exists a Tychonoff star-$\kappa$-compact space which is neither star-$\kappa$-Rothberger nor strongly star-$\kappa$-Menger.
\end{Ex}

\begin{proof}
Since $\kappa$ is weakly compact, $2^\kappa$ is $\kappa$-compact. But $2^\kappa$ is not $\kappa$-Rothberger (see \cite{AC}). Choose $X = 2^\kappa \oplus Y(\kappa^\kappa)$, where $Y(\kappa^\kappa)$ is defined in Example~\ref{ex5}. Then $X$ is a Tychonoff space which is neither star-$\kappa$-Rothberger nor strongly star-$\kappa$-Menger. By Lemma~\ref{lemma8}, $X$ is star-$\kappa$-compact.
\end{proof}

\begin{Ex}
\label{ex9}
Assume that $\kappa$ is a successor cardinal satisfying $\kappa^{<\kappa} = \kappa$ and $\dlk< \dk$. Then there exists a Tychonoff star-$\kappa$-Rothberger space which is neither star-$\kappa$-Scheepers nor strongly star-$\kappa$-Menger.
\end{Ex}

\begin{proof}
Under the assumption that $\kappa$ is a successor cardinal satisfying $\kappa^{<\kappa} = \kappa$ and $\dlk< \dk$, there exists a $\kappa$-Rothberger subspace $Z$ of $\kappa^\kappa$ which is not $\kappa$-Scheepers (see  \cite{AC}).

Choose $X = Y(\kappa^\kappa) \oplus Z$, where $Y(\kappa^\kappa)$ is defined in Example~\ref{ex5}. Then $X$ is a Tychonoff space which is neither star-$\kappa$-Scheepers nor strongly star-$\kappa$-Menger. By Lemma~\ref{lemma10}, $X$ is star-$\kappa$-Rothberger.
\end{proof}

We now consider the reverse implications along the second row of Figure~\ref{fig1}. Recall that a space $X$ is said to be $\omega$-bounded if the closure of every countable subset of $X$ is compact. Note that the following implications hold. \[\omega\text{-bounded} \longrightarrow \text{countably compact} \longrightarrow \text{strongly star-$\kappa$-compact}\]

\begin{Ex}
\label{ex10}
Assume that $\kappa$ is regular. There exists a Tychonoff strongly star-$\kappa$-Hurewicz space which is neither star-$\kappa$-compact nor $\kappa$-Menger.
\end{Ex}

\begin{proof}
Let $Y = \{y_\xi : \xi< \kappa\}$ be a discrete space of cardinality $\kappa$. Then $Y$ is a $\kappa$-Hurewicz space which is not $\kappa$-compact. Consider $Z = [0,\kappa)^\kappa$ endowed with the product topology, where $[0,\kappa)$ carries the order topology. Then $Z$ is $\omega$-bounded.

We claim that $Z$ is not $\kappa$-Menger. For each $\alpha, \xi< \kappa$, put $U_\xi^{(\alpha)} = \{x\in Z : x(\alpha) \leq \xi\}$, and let $\mathcal{U}_\alpha = \{U_\xi^{(\alpha)} : \xi< \kappa\}$. Clearly, $\mathcal{U}_\alpha \in \mathcal{O}_\kappa(Z)$ for every $\alpha< \kappa$. Let $\mathcal{V}_\alpha \in [\mathcal{U}_\alpha]^{<\kappa}$ be arbitrary for every $\alpha< \kappa$. Define $g\in Z$ by
\[g(\alpha) = \begin{cases}
\sup \{\xi+1 : U_\xi^{(\alpha)} \in \mathcal{V}_\alpha\}, & \mbox{if }  \mathcal{V}_\alpha \neq \emptyset\\
0, & \mbox{otherwise}.
\end{cases}\]
If $U_\xi^{(\alpha)} \in \mathcal{V}_\alpha$, then $g(\alpha) > \xi$, and therefore $g\notin U_\xi^{(\alpha)}$. Consequently, $g\notin \bigcup_{\alpha< \kappa} \bigcup \mathcal{V}_\alpha$. Thus $\bigcup_{\alpha< \kappa} \mathcal{V}_\alpha$ does not cover $Z$. Hence $Z$ is not $\kappa$-Menger.

Now choose $X = Y\oplus Z$. Then $X$ is a Tychonoff space which is neither star-$\kappa$-compact and nor $\kappa$-Menger. By Lemma~\ref{lemma4}, $X$ is strongly star-$\kappa$-Hurewicz.
\end{proof}

\begin{Ex}
\label{ex11}
Assume that $\kappa^{<\kappa} = \kappa$ and $\bk< \dk$. Then there exists a Tychonoff strongly star-$\kappa$-Scheepers space which is neither star-$\kappa$-Hurewicz nor $\kappa$-Menger.
\end{Ex}

\begin{proof}
Under the assumption $\kappa^{<\kappa} = \kappa$ and $\bk< \dk$, there exists a $\kappa$-Scheepers subspace $Y$ of $\kappa^\kappa$ which is not $\kappa$-Hurewicz (see \cite{AC}). Let $Z = [0,\kappa)^\kappa$ be endowed with the product topology, where $[0,\kappa)$ carries the order topology. Then $Z$ is $\omega$-bounded but not $\kappa$-Menger (see Example~\ref{ex10}).

Put $X = Y\oplus Z$. Clearly, $X$ is a Tychonoff space which is not star-$\kappa$-Hurewicz. By Lemma~\ref{lemma6}, $X$ is strongly star-$\kappa$-Scheepers.
\end{proof}

\begin{Ex}
\label{ex14}
Assume that $\kappa^{<\kappa} = \kappa$ and $\dlk< \dk$. Then there exists a Tychonoff strongly star-$\kappa$-Menger space which is neither star-$\kappa$-Scheepers nor $\kappa$-Menger.
\end{Ex}

\begin{proof}
Under the assumption $\kappa^{<\kappa} = \kappa$ and $\dlk< \dk$, there exists a $\kappa$-Menger subspace $Y$ of $\kappa^\kappa$ which is not $\kappa$-Scheepers (see \cite{AC}). Let $Z = [0,\kappa)^\kappa$ be endowed with the product topology, where $[0,\kappa)$ carries the order topology. Then $Z$ is $\omega$-bounded but not $\kappa$-Menger (see Example~\ref{ex10}).

Put $X = Y\oplus Z$. Clearly, $X$ is a Tychonoff space which is neither star-$\kappa$-Scheepers nor $\kappa$-Menger. By Lemma~\ref{lemma7}, $X$ is strongly star-$\kappa$-Menger.
\end{proof}

\begin{Ex}
\label{ex13}
Assume that $\kappa$ is a weakly compact cardinal. Then there exists a Tychonoff strongly star-$\kappa$-compact space which is neither star-$\kappa$-Rothberger nor $\kappa$-Menger.
\end{Ex}

\begin{proof}
Recall from \cite{AC} that the $\kappa$-Cantor space $2^\kappa$ is not $\kappa$-Rothberger. Since $\kappa$ is weakly compact, $2^\kappa$ is $\kappa$-compact.  Let $Y = [0,\kappa)^\kappa$ be endowed with the product topology, where $[0,\kappa)$ carries the order topology. Then $Y$ is $\omega$-bounded but not $\kappa$-Menger (see Example~\ref{ex10}).

Put $X = 2^\kappa \oplus Y$. Clearly, $X$ is a Tychonoff space which is neither star-$\kappa$-Rothberger nor $\kappa$-Menger. By Lemma~\ref{lemma8}, $X$ is strongly star-$\kappa$-compact.
\end{proof}

\begin{Ex}
\label{ex12}
Assume that $\kappa$ is a successor cardinal satisfying $\kappa^{<\kappa} = \kappa$ and $\dlk< \dk$. Then there exists a $T_1$ strongly star-$\kappa$-Rothberger space which is neither star-$\kappa$-Scheepers nor $\kappa$-Rothberger.
\end{Ex}

\begin{proof}
Under the assumption that $\kappa$ is a successor cardinal satisfying $\kappa^{<\kappa} = \kappa$ and $\dlk< \dk$, there exists a $\kappa$-Rothberger subspace $Y$ of $\kappa^\kappa$ which is not $\kappa$-Scheepers (see  \cite{AC}). Let $Z$ be the space considered in Example~\ref{ex3}. Then $Z$ is a $T_1$ strongly star-$\kappa$-Rothberger space which is not $\kappa$-Rothberger.

Put $X = Y\oplus Z$. Clearly, $X$ is a $T_1$ space which is neither star-$\kappa$-Scheepers nor $\kappa$-Rothberger. By Lemma~\ref{lemma10}, $X$ is strongly star-$\kappa$-Rothberger.
\end{proof}

Finally, we consider the reverse implications along the columns of Figure~\ref{fig1}.

\begin{Ex}
\label{ex1}
For every cardinal $\kappa$, there exists a $T_1$ strongly star-$\kappa$-compact space which is not $\kappa$-Menger.
\end{Ex}

\begin{proof}
Define a topology on $X = \omega \cup \kappa^\kappa$ by declaring every point of $\omega$ isolated, and for each $f\in \kappa^\kappa$, taking $\{\{f\} \cup (\omega \setminus F) : F\in [\omega]^{<\omega}\}$ as a neighbourhood base at $f$. Clearly, $X$ is a $T_1$ space. Moreover, $\omega$ is dense in $X$. Then $X$ is strongly star-$\kappa$-compact.

We now show that $X$ is not $\kappa$-Menger. For each $\alpha, \xi< \kappa$, let \[U^{(\alpha)}_\xi = \omega \cup \{f\in \kappa^\kappa : f(\alpha) = \xi\},\] and put $\mathcal{U}_\alpha = \{U^{(\alpha)}_\xi : \xi< \kappa\}$. Each $U^{(\alpha)}_\xi$ is open in $X$. Indeed, every point of $\omega$ is isolated, and if $f\in \kappa^\kappa$ satisfies $f(\alpha) = \xi$, then $\{f\} \cup \omega \subseteq U^{(\alpha)}_\xi$. Moreover, $|\mathcal{U}_\alpha| = \kappa$ and $\mathcal{U}_\alpha$ covers $X$. Hence $\mathcal{U}_\alpha \in \mathcal{O}_\kappa(X)$ for every $\alpha< \kappa$.

Suppose that for each $\alpha< \kappa$, $\mathcal{V}_\alpha \in [\mathcal{U}_\alpha]^{<\kappa}$. Let $A_\alpha = \{\xi< \kappa : U^{(\alpha)}_\xi \in \mathcal{V}_\alpha\}$. Then $|A_\alpha| < \kappa$, and hence we may choose $g(\alpha) \in \kappa \setminus A_\alpha$ for each $\alpha< \kappa$. Thus $g\in \kappa^\kappa$. For each $\alpha< \kappa$ and each $U^{(\alpha)}_\xi \in \mathcal{V}_\alpha$, we have $\xi\in A_\alpha$, and therefore $g(\alpha) \neq \xi$. Consequently, $g\notin U^{(\alpha)}_\xi$. It follows that $g\notin \bigcup_{\alpha< \kappa} \bigcup \mathcal{V}_\alpha$. Hence $X\neq \bigcup_{\alpha< \kappa} \bigcup \mathcal{V}_\alpha$. Therefore, $X$ is not $\kappa$-Menger.
\end{proof}

\begin{Ex}
\label{ex3}
For every cardinal $\kappa$, there exists a $T_1$ strongly star-$\kappa$-Rothberger space which is not $\kappa$-Rothberger.
\end{Ex}

\begin{proof}
Let $\lambda = (2^\kappa)^+$, and let $X$ be a set of cardinality $\lambda$ such that $\kappa^\kappa \subseteq X$. Endow $X$ with the co-<$\lambda$ topology, that is, \[\tau = \{\emptyset\} \cup \{U\subseteq X : |X\setminus U|< \lambda\}.\] Clearly, $X$ is a $T_1$ space. We first show that $X$ is strongly star-$\kappa$-Rothberger. Let $(\mathcal{U}_\alpha : \alpha< \kappa) \in \mathcal{O}_\kappa(X)^\kappa$. For each $\alpha< \kappa$, write $\mathcal{U}_\alpha = \{U^{(\alpha)}_\xi : \xi< \kappa\}$. Since $|X\setminus U^{(\alpha)}_\xi| < \lambda$, we have $|X\setminus U^{(\alpha)}_\xi| \leq 2^\kappa$. Hence \[\left|\bigcup_{\xi< \kappa} (X\setminus U^{(\alpha)}_\xi)\right| \leq \kappa \cdot 2^\kappa = 2^\kappa< \lambda.\] Thus $\bigcap_{\xi< \kappa} U^{(\alpha)}_\xi \neq \emptyset$. Choose $x_\alpha \in \bigcap_{\xi< \kappa} U^{(\alpha)}_\xi$. Then every member of $\mathcal{U}_\alpha$ contains $x_\alpha$, and therefore $St(x_\alpha, \mathcal{U}_\alpha) = \bigcup \mathcal{U}_\alpha = X$. Consequently, $X = \bigcup_{\alpha< \kappa} St(x_\alpha, \mathcal{U}_\alpha)$, and hence $X$ is strongly star-$\kappa$-Rothberger.

We now show that $X$ is not $\kappa$-Rothberger. For each $\alpha, \xi< \kappa$, put \[V^{(\alpha)}_\xi = (X\setminus \kappa^\kappa) \cup \{f\in \kappa^\kappa : f(\alpha) = \xi\},\] and let $\mathcal{V}_\alpha = \{V^{(\alpha)}_\xi : \xi< \kappa\}$. Since \[X\setminus V^{(\alpha)}_\xi = \{f\in \kappa^\kappa : f(\alpha) \neq \xi\}\] has cardinality at most $|\kappa^\kappa| = 2^\kappa< \lambda$, each $V^{(\alpha)}_\xi$ is open in $X$. Moreover, $|\mathcal{V}_\alpha| = \kappa$ and $\bigcup \mathcal{V}_\alpha = X$. Thus $\mathcal{V}_\alpha \in \mathcal{O}_\kappa(X)$ for every $\alpha< \kappa$.

Now choose arbitrarily $V^{(\alpha)}_{\xi_\alpha} \in \mathcal{V}_\alpha$, $\alpha< \kappa$. For every $\alpha< \kappa$, choose $g(\alpha) \in\kappa \setminus \{\xi_\alpha\}$. Then $g\in \kappa^\kappa$, and $g(\alpha) \neq \xi_\alpha$ for every $\alpha< \kappa$. Hence $g\notin V^{(\alpha)}_{\xi_\alpha}$ for every $\alpha< \kappa$. Therefore, $g\notin \bigcup_{\alpha< \kappa} V^{(\alpha)}_{\xi_\alpha}$. Thus no choice of one member from each $\mathcal{V}_\alpha$ covers $X$, and consequently $X$ is not $\kappa$-Rothberger.
\end{proof}

\begin{Ex}
\label{ex4}
For every cardinal $\kappa$, there exists a Tychonoff space which is both star-$\kappa$-compact and star-$\kappa$-Rothberger, but not strongly star-$\kappa$-Menger.
\end{Ex}

\begin{proof}
Let $Y(\kappa^\kappa)$ be the space defined in Example~\ref{ex5}. Then $Y(\kappa^\kappa)$ is a Tychonoff star-$\kappa$-compact space which is not strongly star-$\kappa$-Menger. By \cite[Example~3.7]{BCKM}, $X$ is star-$\kappa$-Rothberger.
\end{proof}

\begin{figure}[h]
\begin{adjustbox}{max width=\textwidth,max
height=\textheight,keepaspectratio,center}
\begin{tikzcd}[column sep=7.5ex, row sep=13ex, arrows={crossing over}]
\text{star-$\kappa$-compact} \arrow[rrrr, bend left=12, "/" {marking, scale=0.6, pos=0.7}, "\text{Ex.~\ref{ex8}}"] \arrow[r, shift left] \arrow[d, shift left, "/" {marking, scale=0.6, pos=0.51}, "\text{Ex.~\ref{ex4}}", sloped, pos=0.55] & \text{star-$\kappa$-Hurewicz} \arrow[l, shift left, "/" {marking, scale=0.6}, "\text{Ex.~\ref{ex5}}"] \arrow[r, shift left] \arrow[d, shift left, "/" {marking, scale=0.6, pos=0.41}, "\text{Ex.~\ref{ex5}}", sloped, pos=0.45] & \text{star-$\kappa$-Scheepers} \arrow[l, shift left, "/" {marking, scale=0.6}, "\text{Ex.~\ref{ex6}}"] \arrow[r, shift left] \arrow[d, shift left, "/" {marking, scale=0.6, pos=0.8}, "\text{Ex.~\ref{ex5}}", sloped, pos=0.3] & \text{star-$\kappa$-Menger} \arrow[l, shift left, "/" {marking, scale=0.6}, "\text{Ex.~\ref{ex7}}"] \arrow[r, shift right, "/" {marking, scale=0.6}, "\text{Ex.~\ref{ex8}}", swap] \arrow[d, shift left, "/" {marking, scale=0.6, pos=0.51}, "\text{Ex.~\ref{ex5}}", sloped, pos=0.55] & \text{star-$\kappa$-Rothberger} \arrow[ll, bend left=15, "/" {marking, pos=0.7, scale=0.6}, "\text{Ex.~\ref{ex9}}", pos=0.7, sloped, swap] \arrow[l, shift right] \arrow[d, shift left, "/" {marking, scale=0.6, pos=0.51}, "\text{Ex.~\ref{ex4}}", sloped, pos=0.55] \\
\text{\parbox{3.5cm}{\centering strongly\\ star-$\kappa$-compact}} \arrow[rrrr, bend left=12, "/" {marking, scale=0.6, pos=0.7}, "\text{Ex.~\ref{ex13}}", pos=0.65, sloped] \arrow[r, shift left] \arrow[u, shift left] \arrow[d, shift left, "/" {marking, scale=0.6, pos=0.51}, "\text{Ex.~\ref{ex1}}", sloped, pos=0.55] & \text{\parbox{3.5cm}{\centering strongly\\ star-$\kappa$-Hurewicz}} \arrow[l, shift left, "/" {marking, scale=0.6}, "\text{Ex.~\ref{ex10}}"] \arrow[r, shift left] \arrow[u, shift left] \arrow[d, shift left, "/" {marking, scale=0.6, pos=0.35}, "\text{Ex.~\ref{ex10}}", sloped, pos=0.39] & \text{\parbox{3.5cm}{\centering strongly\\ star-$\kappa$-Scheepers}} \arrow[l, shift left, "/" {marking, scale=0.6}, "\text{Ex.~\ref{ex11}}"] \arrow[r, shift left] \arrow[u, shift left] \arrow[d, shift left, "/" {marking, scale=0.6, pos=0.8}, "\text{Ex.~\ref{ex10}}", sloped, pos=0.25] & \text{\parbox{3.5cm}{\centering strongly\\ star-$\kappa$-Menger}} \arrow[l, shift left, "/" {marking, scale=0.6}, "\text{Ex.~\ref{ex14}}"] \arrow[r, shift right, "/" {marking, scale=0.6}, "\text{Ex.~\ref{ex13}}", swap] \arrow[u, shift left] \arrow[d, shift left, "/" {marking, scale=0.6, pos=0.42}, "\text{Ex.~\ref{ex10}}", sloped, pos=0.46] & \text{\parbox{3.5cm}{\centering strongly\\ star-$\kappa$-Rothberger}} \arrow[ll, bend left=15, "/" {marking, pos=0.7, scale=0.6}, "\text{Ex.~\ref{ex12}}", pos=0.735, sloped, swap] \arrow[l, shift right] \arrow[u, shift left] \arrow[d, shift left, "/" {marking, scale=0.6, pos=0.51}, "\text{Ex.~\ref{ex3}}", sloped, pos=0.55] \\
\text{$\kappa$-compact} \arrow[rrrr, bend left=12, "/" {marking, scale=0.6, pos=0.7}, "\text{\cite{AC}}", pos=0.65, sloped] \arrow[r, shift left] \arrow[u, shift left] & \text{$\kappa$-Hurewicz} \arrow[l, shift left, "/" {marking, scale=0.6,pos=0.5}, "\text{\cite{AC}}", pos=0.65] \arrow[r,"\text{$\kappa$ is regular}", shift left] \arrow[u, shift left] & \text{$\kappa$-Scheepers} \arrow[l, shift left, "/" {marking, scale=0.6, pos=0.5}, "\text{\cite{AC}}", pos=0.65] \arrow[r, shift left] \arrow[u, shift left] & \text{$\kappa$-Menger} \arrow[l, shift left, "/" {marking, scale=0.6, pos=0.5}, "\text{\cite{AC}}", pos=0.4] \arrow[r, shift right, "/" {marking, scale=0.6, pos=0.5}, "\text{\cite{AC}}", swap, pos=0.35] \arrow[u, shift left] & \text{$\kappa$-Rothberger} \arrow[ll, bend left=15, "/" {marking, pos=0.7, scale=0.6}, "\text{\cite{AC}}"] \arrow[l, shift right] \arrow[u, shift left]
\end{tikzcd}
\end{adjustbox}
\caption{Implications and non-implications among $\kappa$-covering properties and their star variants.}
\label{fig2}
\end{figure}
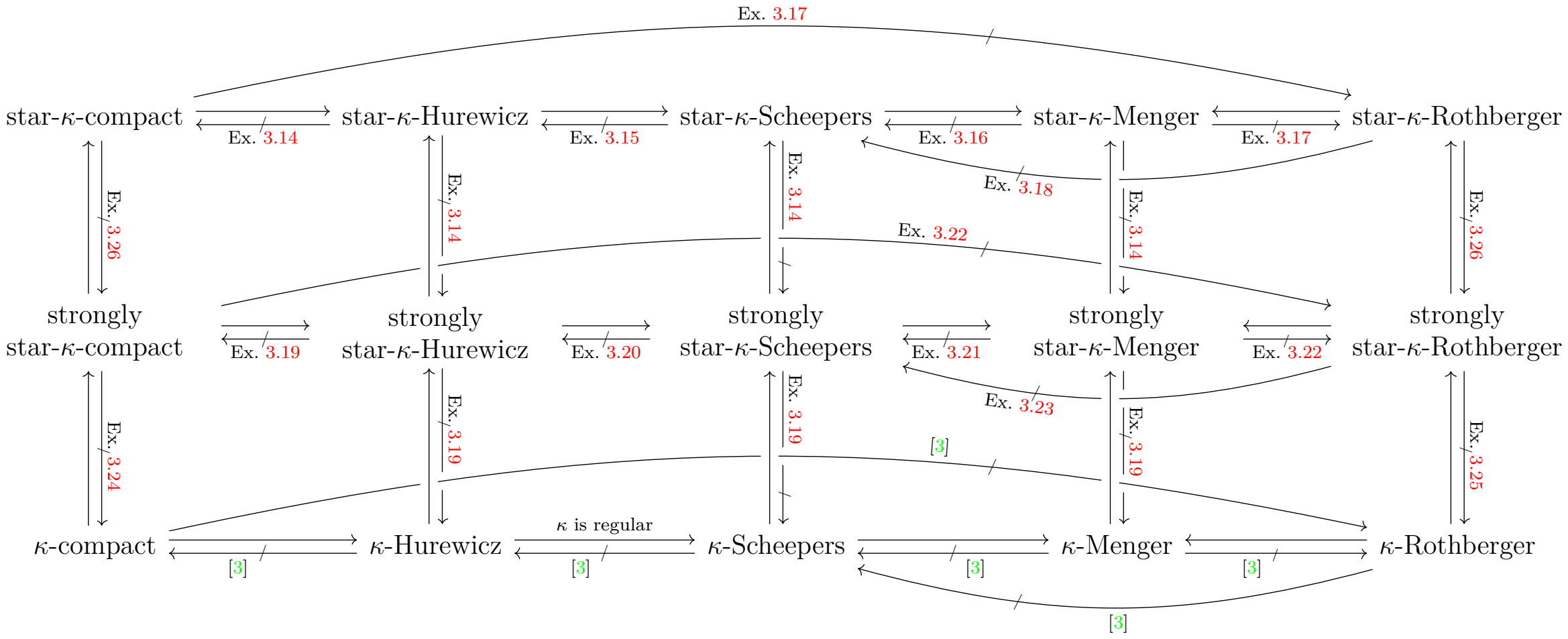

For the convenience of the reader, we now reproduce Figure~\ref{fig1} as Figure~\ref{fig2}, incorporating the non-implications discussed in this section.

\section{Combinatorial aspects of uncountable star selection principles in generalized $\Psi$-spaces}
\label{sca}
Recall that a family $\mathcal{A} \subseteq [\kappa]^\kappa$ is called an almost disjoint family if $|A\cap B|< \kappa$ for any distinct $A,B\in \mathcal{A}$. It is called maximal almost disjoint (MAD) family if for every $B\in [\kappa]^\kappa \setminus \mathcal{A}$, there exists $A\in \mathcal{A}$ such that $|A\cap B| = \kappa$, in other words, if $\mathcal{A}$ is not included in any larger almost disjoint family. Following \cite{xss20}, the generalized Isbell--Mr\'{o}wka space (or generalized $\Psi$-space) $\Psi_\kappa(\mathcal{A}) = \kappa \cup \mathcal{A}$ associated with an almost disjoint family $\mathcal{A}$ is topologized as follows. Each point of $\kappa$ is isolated, and a neighbourhood base at each $A\in \mathcal{A}$ consists of sets of the form $\{A\} \cup (A\setminus F)$, where $F\in [A]^{<\kappa}$. Thus a neighbourhood base at $A\in \mathcal{A}$ is given by $\mathcal{B}(A) = \{\{A\} \cup (A\setminus F) : F\in [A]^{<\kappa}\}$. This construction is the case $\lambda = \kappa$ of the generalized Mr\'{o}wka spaces considered in \cite{xss20}. If $\kappa$ is regular, then the sets $N(A,\beta) = \{A\} \cup (A\setminus \beta)$, $\beta< \kappa$, also form a neighbourhood base at $A$. Clearly, $\Psi_\kappa(\mathcal{A})$ is a zero-dimensional Hausdorff space, and therefore is a Tychonoff space.

In the countable setting, star selection principles in classical $\Psi$-spaces \cite{mrowka54} have been extensively studied, both from a combinatorial point of view and in connection with the construction of counterexamples (see \cite{tsaban15, bm09, RBS, ca23, ca24}). In view of these developments, having introduced star selection principles in the uncountable setting, we naturally turn our attention in this section to generalized $\Psi$-spaces as a first step towards understanding the behaviour of these principles.

\subsection{Combinatorial characterizations}
\label{scc}
In this section, we establish combinatorial characterizations of star-$\kappa$-Menger, star-$\kappa$-Scheepers, star-$\kappa$-Hurewicz, and star-$\kappa$-Rothberger $\Psi_\kappa$-spaces (Theorems~\labelcref{thm1,thm4,thm2,thm3}), generalizing \cite[Theorem~2.1]{tsaban15}, \cite[Theorem~3.12]{ca23}, \cite[Theorem~2.2]{tsaban15}, and \cite[Theorem~4.3]{tsaban15}, respectively, to the uncountable setting. We present the proof only for the star-$\kappa$-Scheepers characterization; the proofs of the other theorems are analogous and hence omitted.

Let $\mathcal{A} \subseteq [\kappa]^\kappa$ be an almost disjoint family. For an assignment $A\mapsto f_A$ from $\mathcal{A}$ into $\kappa^\kappa$, and for each $\alpha< \kappa$, a partition $\mathcal{P}_\alpha$ of $\mathcal{A}$ with $|\mathcal{P}_\alpha| \leq \kappa$, put $T_\alpha(P) = \bigcup_{B\in P} (B\setminus f_B(\alpha))$, $P\in \mathcal{P}_\alpha$. For $A\in \mathcal{A}$, let $P_\alpha(A)$ denote the unique member of $\mathcal{P}_\alpha$ containing $A$.

\begin{Th}
\label{thm1}
Assume that $\kappa$ is regular. Let $\mathcal{A} \subseteq [\kappa]^\kappa$ be an almost disjoint family. Then the following assertions are equivalent.
\begin{enumerate}[wide=0pt,label={\upshape(\arabic*)},leftmargin=*]
\item $\Psi_\kappa(\mathcal{A})$ is star-$\kappa$-Menger.
\item For every assignment $A\mapsto f_A$ from $\mathcal{A}$ into $\kappa^\kappa$ and every sequence $(\mathcal{P}_\alpha : \alpha< \kappa)$ of partitions of $\mathcal{A}$ with $|\mathcal{P}_\alpha| \leq \kappa$, there are $\mathcal{F}_\alpha \in [\mathcal{P}_\alpha]^{<\kappa}$, $\alpha< \kappa$, such that for each $A\in \mathcal{A}$, there exists $\alpha< \kappa$ for which \[T_\alpha(P_\alpha(A)) \cap \bigcup_{P\in \mathcal{F}_\alpha} T_\alpha(P) \neq \emptyset.\]
\end{enumerate}
\end{Th}

To establish the desired characterization of star-$\kappa$-Scheepers $\Psi_\kappa$-spaces, we first prove the following useful lemma.

\begin{Lemma}
\label{lemma23}
Let $\mathcal{A} \subseteq [\kappa]^\kappa$ be an almost disjoint family. The following assertions are equivalent.
\begin{enumerate}[wide=0pt,label={\upshape(\arabic*)},leftmargin=*]
\item For every assignment $A\mapsto f_A$ from $\mathcal{A}$ into $\kappa^\kappa$ and  every sequence $(\mathcal{P}_\alpha : \alpha< \kappa)$ of partitions of $\mathcal{A}$ with $|\mathcal{P}_\alpha| \leq \kappa$, there are $\mathcal{F}_\alpha \in [\mathcal{P}_\alpha]^{<\kappa}$, $\alpha< \kappa$, such that for every $\mathcal{F} \in [\mathcal{A}]^{<\kappa}$, there exists $\alpha< \kappa$ such that \[T_\alpha(P_\alpha(A)) \cap \bigcup_{P\in \mathcal{F}_\alpha} T_\alpha(P) \neq \emptyset\] for every $A\in \mathcal{F}$.
\item For every assignment $A\mapsto f_A$ from $\mathcal{A}$ into $\kappa^\kappa$ and  every sequence $(\mathcal{P}_\alpha : \alpha< \kappa)$ of partitions of $\mathcal{A}$ with $|\mathcal{P}_\alpha| \leq \kappa$, there are $\mathcal{F}_\alpha \in [\mathcal{P}_\alpha]^{<\kappa}$, $\alpha< \kappa$, such that for every $\mathcal{F} \in [\mathcal{A}]^{<\kappa}$, the size of \[\left\{\alpha< \kappa : (\forall A\in \mathcal{F})\, T_\alpha(P_\alpha(A)) \cap \bigcup_{P\in \mathcal{F}_\alpha} T_\alpha(P) \neq \emptyset\right\}\] is $\kappa$.
\end{enumerate}
\end{Lemma}

\begin{proof}
$(1)\Rightarrow (2)$. Fix an assignment $A\mapsto f_A$ from $\mathcal{A}$ into $\kappa^\kappa$ and a sequence $(\mathcal{P}_\alpha : \alpha< \kappa)$ of partitions of $\mathcal{A}$ with $|\mathcal{P}_\alpha| \leq \kappa$. Partition $\kappa$ into pairwise disjoint sets $(I_\xi : \xi< \kappa)$ of cardinality $\kappa$. For each $\xi< \kappa$, let $e_\xi: \kappa\to I_\xi$ be the increasing
enumeration of $I_\xi$. Define $f_A^{(\xi)}(\beta) = f_A(e_\xi(\beta))$ and $\mathcal{P}_\beta^{(\xi)} = \mathcal{P}_{e_\xi(\beta)}$ for every $A\in \mathcal{A}$ and $\beta< \kappa$. By $(1)$, applied to the assignment $A\mapsto f_A^{(\xi)}$ and the sequence $(\mathcal{P}_\beta^{(\xi)} : \beta< \kappa)$, there exist $\mathcal{F}_\beta^{(\xi)} \in [\mathcal{P}_\beta^{(\xi)}]^{<\kappa}$, $\beta< \kappa$, such that for every $\mathcal{F}\in [\mathcal{A}]^{<\kappa}$, there is some $\beta< \kappa$ for which \[T_{e_\xi(\beta)}(P_{e_\xi(\beta)}(A)) \cap \bigcup_{P\in \mathcal{F}_\beta^{(\xi)}} T_{e_\xi(\beta)}(P) \neq \emptyset\] for every $A\in \mathcal{F}$.

Since the sets $I_\xi$ partition $\kappa$, we may define $\mathcal{F}_\alpha = \mathcal F_\beta^{(\xi)}$ whenever $\alpha = e_\xi(\beta)$. Then $\mathcal{F}_\alpha \in [\mathcal{P}_\alpha]^{<\kappa}$ for every $\alpha< \kappa$. Now fix $\mathcal{F}\in [\mathcal{A}]^{<\kappa}$ and let \[G_{\mathcal{F}} = \left\{\alpha< \kappa : (\forall A\in \mathcal{F})\, T_\alpha(P_\alpha(A)) \cap \bigcup_{P\in \mathcal{F}_\alpha} T_\alpha(P) \neq \emptyset\right\}.\] For every $\xi< \kappa$, the preceding application of $(1)$ gives an index $\alpha_\xi \in I_\xi \cap G_{\mathcal{F}}$. Since the sets $I_\xi$ are pairwise disjoint, the indices $(\alpha_\xi : \xi< \kappa)$ are distinct. Consequently, $|G_{\mathcal{F}}| = \kappa$.
\end{proof}

\begin{Th}
\label{thm4}
Assume that $\kappa$ is regular. Let $\mathcal{A} \subseteq [\kappa]^\kappa$ be an almost disjoint family. The following assertions are equivalent.
\begin{enumerate}[wide=0pt,label={\upshape(\arabic*)},leftmargin=*]
\item $\Psi_\kappa(\mathcal{A})$ is star-$\kappa$-Scheepers.
\item For every assignment $A\mapsto f_A$ from $\mathcal{A}$ into $\kappa^\kappa$ and  every sequence $(\mathcal{P}_\alpha : \alpha< \kappa)$ of partitions of $\mathcal{A}$ with $|\mathcal{P}_\alpha| \leq \kappa$, there are $\mathcal{F}_\alpha \in [\mathcal{P}_\alpha]^{<\kappa}$, $\alpha< \kappa$, such that for every $\mathcal{F} \in [\mathcal A]^{<\kappa}$, there exists $\alpha< \kappa$ such that \[T_\alpha(P_\alpha(A)) \cap \bigcup_{P\in \mathcal{F}_\alpha} T_\alpha(P) \neq \emptyset\] for every $A\in \mathcal{F}$.
\end{enumerate}
\end{Th}

\begin{proof}
Assume first that $\Psi_\kappa(\mathcal{A})$ is star-$\kappa$-Scheepers. Fix an assignment $A\mapsto f_A$ from $\mathcal{A}$ into $\kappa^\kappa$ and a sequence
$(\mathcal{P}_\alpha : \alpha< \kappa)$ as in $(2)$. For $P\in \mathcal{P}_\alpha$, define $W_\alpha(P) = P\cup T_\alpha(P)$. Clearly, $W_\alpha(P)$ is open in $\Psi_\kappa(\mathcal{A})$. Now put \[\mathcal{U}_\alpha = \{W_\alpha(P) : P\in \mathcal{P}_\alpha\} \cup \{\{\xi\} : \xi< \kappa\}.\] Then $\mathcal{U}_\alpha \in \mathcal{O}_\kappa(\Psi_\kappa(\mathcal{A}))$. Choose $\mathcal{V}_\alpha \in [\mathcal{U}_\alpha]^{<\kappa}$ witnessing the star-$\kappa$-Scheepers property. We may replace every selected singleton $\{\xi\} \in \mathcal{V}_\alpha$ for which $\xi\in T_\alpha(P)$ for some $P\in \mathcal{P}_\alpha$ by one such $W_\alpha(P)$. If no such $P$ exists, then $St(\{\xi\}, \mathcal{U}_\alpha)$ contains no point of $\mathcal{A}$, and hence the singleton may simply be discarded when considering the coverage of $\mathcal{A}$.

Accordingly, there are $\mathcal{F}_\alpha \in [\mathcal{P}_\alpha]^{<\kappa}$, $\alpha< \kappa$, such that $\bigcup \mathcal{V}_\alpha = \bigcup_{P\in \mathcal{F}_\alpha} W_\alpha(P)$ as far as the coverage of $\mathcal{A}$ is concerned. Let $\mathcal{F}\in [\mathcal{A}]^{<\kappa}$. Choose $\alpha< \kappa$ such that \[\mathcal{F} \subseteq St\left(\bigcup_{P\in \mathcal{F}_\alpha} W_\alpha(P), \mathcal{U}_\alpha\right).\] The only member of $\mathcal{U}_\alpha$ containing the non-isolated point $A\in \mathcal{F}$ is $W_\alpha(P_\alpha(A))$. Hence \[W_\alpha(P_\alpha(A)) \cap \bigcup_{P\in \mathcal{F}_\alpha} W_\alpha(P) \neq \emptyset\] for every $A\in \mathcal{F}$. Since the members of $\mathcal{P}_\alpha$ are pairwise disjoint, an intersection between two distinct $W_\alpha$-sets can occur only in their isolated parts. It follows that \[T_\alpha(P_\alpha(A)) \cap \bigcup_{P\in\mathcal F_\alpha}T_\alpha(P) \neq \emptyset\] for every $A\in \mathcal{F}$.

Conversely, suppose that (2) holds. Let $(\mathcal{U}_\alpha : \alpha< \kappa) \in \mathcal{O}_\kappa(\Psi_\kappa(\mathcal{A}))^\kappa$. Write $\mathcal{U}_\alpha = \{U_\xi^{(\alpha)} : \xi< \kappa\}$. For every $A\in \mathcal{A}$ and every $\alpha< \kappa$, choose $\gamma_\alpha(A) < \kappa$ and $f_A(\alpha) < \kappa$ so that $N(A,f_A(\alpha)) \subseteq U_{\gamma_\alpha(A)}^{(\alpha)}$. Let $\mathcal{P}_\alpha$ be the partition of $\mathcal{A}$ consisting of the nonempty fibers of $\gamma_\alpha$. Thus if
$P\in \mathcal{P}_\alpha$, then for some $\xi< \kappa$, $P = \{A\in \mathcal{A} : \gamma_\alpha(A) = \xi\}$. In particular, $P\cup T_\alpha(P) \subseteq U_\xi^{(\alpha)}$.

Apply (2), and choose $\mathcal{F}_\alpha \in [\mathcal{P}_\alpha]^{<\kappa}$. For each $P\in \mathcal{F}_\alpha$, let $\xi(P)< \kappa$ be such that $P = \{A\in \mathcal{A} : \gamma_\alpha(A) = \xi(P)\}$, and put $\mathcal{H}_\alpha = \{U_{\xi(P)}^{(\alpha)} : P\in \mathcal{F}_\alpha\}$. We also take care of the isolated part $\kappa$. For every $\xi\leq \alpha$, choose $U(\alpha,\xi)\in \mathcal{U}_\alpha$ such that $\xi\in U(\alpha,\xi)$, and put \[\mathcal{V}_\alpha = \mathcal{H}_\alpha \cup \{U(\alpha,\xi) : \xi\leq \alpha\}.\] Clearly, $|\mathcal{V}_\alpha| < \kappa$. Let $\mathcal{F}\in [\Psi_\kappa(\mathcal{A})]^{<\kappa}$. By Lemma~\ref{lemma23}, the size of \[B = \left\{\alpha< \kappa : (\forall A\in \mathcal{F} \cap \mathcal{A})\, T_\alpha(P_\alpha(A)) \cap \bigcup_{P\in \mathcal{F}_\alpha} T_\alpha(P) \neq \emptyset\right\}\] is $\kappa$. Since $P_\alpha(A) \cup T_\alpha(P_\alpha(A)) \subseteq U_{\gamma_\alpha(A)}^{(\alpha)}$ for every $A\in \mathcal{F} \cap \mathcal{A}$ and $\bigcup_{P\in \mathcal{F}_\alpha} T_\alpha(P)$ is contained in $\bigcup \mathcal{V}_\alpha$, we obtain $\mathcal{F} \cap \mathcal{A} \subseteq St(\bigcup \mathcal{V}_\alpha, \mathcal{U}_\alpha)$ for every $\alpha\in B$.

Since $\kappa$ is regular and $|\mathcal{F} \cap \kappa|< \kappa$, there is
$\gamma\in B$ such that $\mathcal{F} \cap \kappa \subseteq \gamma$. It follows that $\mathcal{F} \cap \kappa \subseteq \bigcup \mathcal{V}_\gamma \subseteq St(\bigcup \mathcal{V}_\gamma, \mathcal{U}_\gamma)$. Thus there exists $\alpha< \kappa$ such that $\mathcal{F} \subseteq St(\bigcup \mathcal{V}_\alpha, \mathcal{U}_\alpha)$, and hence $\Psi_\kappa(\mathcal{A})$ is star-$\kappa$-Scheepers.
\end{proof}

\begin{Th}
\label{thm2}
Assume that $\kappa$ is regular. Let $\mathcal{A} \subseteq [\kappa]^\kappa$ be an almost disjoint family. The following assertions are equivalent.
\begin{enumerate}[wide=0pt,label={\upshape(\arabic*)},leftmargin=*]
\item $\Psi_\kappa(\mathcal{A})$ is star-$\kappa$-Hurewicz.
\item For every assignment $A\mapsto f_A$ from $\mathcal{A}$ into $\kappa^\kappa$ and  every sequence $(\mathcal{P}_\alpha : \alpha< \kappa)$ of partitions of $\mathcal{A}$ with $|\mathcal{P}_\alpha| \leq \kappa$, there are $\mathcal{F}_\alpha \in [\mathcal{P}_\alpha]^{<\kappa}$, $\alpha< \kappa$, such that for every $A\in \mathcal{A}$, the size of \[\left\{\alpha< \kappa : T_\alpha(P_\alpha(A)) \cap \bigcup_{P\in \mathcal{F}_\alpha} T_\alpha(P) = \emptyset\right\}\] is less than $\kappa$.
\end{enumerate}
\end{Th}

\begin{Th}
\label{thm3}
Assume that $\kappa$ is regular. Let $\mathcal{A} \subseteq [\kappa]^\kappa$ be an almost disjoint family. The following assertions are equivalent.
\begin{enumerate}[wide=0pt,label={\upshape(\arabic*)},leftmargin=*]
\item $\Psi_\kappa(\mathcal{A})$ is star-$\kappa$-Rothberger.
\item For every assignment $A\mapsto f_A$ from $\mathcal{A}$ into $\kappa^\kappa$ and every sequence $(\mathcal{P}_\alpha : \alpha< \kappa)$ of partitions of $\mathcal{A}$ with $|\mathcal{P}_\alpha| \leq \kappa$, there are $Q_\alpha \in \mathcal{P}_\alpha$, $\alpha< \kappa$, such that for every $A\in \mathcal{A}$, there exists $\alpha< \kappa$ for which $T_\alpha(P_\alpha(A)) \cap T_\alpha(Q_\alpha) \neq \emptyset$.
\end{enumerate}
\end{Th}

\subsection{Cardinal characterizations}
\label{scc1}
Unlike in the classical countable setting, we establish that for generalized $\Psi$-spaces, the $\kappa$-Menger, $\kappa$-Scheepers, and $\kappa$-Hurewicz properties are each equivalent to their corresponding strongly star variants and to cardinality bounds involving $\dk$, $\dlk$, and $\bk$, respectively. The following lemma will be useful in establishing these results.

\begin{Lemma}
\label{lemma18}
Assume that $\kappa$ is regular. If $|X|< \dk$ (respectively, $\dlk$, $\bk$), then $X$ is $\kappa$-Menger (respectively, $\kappa$-Scheepers, $\kappa$-Hurewicz).
\end{Lemma}

\begin{proof}
We prove only the $\kappa$-Scheepers case, as the remaining cases follow by analogous arguments. Let $(\mathcal{U}_\alpha : \alpha< \kappa) \in \mathcal{O}_\kappa(X)^\kappa$ be such that no $\mathcal{U}_\alpha$ has a subcover of cardinality less than $\kappa$. Partition $\kappa$ into two disjoint sets $A$ and $B$ such that $|A| = |B| = \kappa$, and fix a bijection $e: \kappa\to A$. Write $\mathcal{U}_{e(\xi)} = \{U^{(e(\xi))}_\eta : \eta< \kappa\}$ for each $\xi< \kappa$. For every $x\in X$, define $f_x\in \kappa^\kappa$ by \[f_x(\xi) = \min\{\eta< \kappa : x\in U^{(e(\xi))}_\eta\}.\] Put $Y = \{f_x : x\in X\}$. Since $|Y|< \dlk$, the family $\maxk(Y)$ is not dominating (see Lemma~\ref{lemma5}), i.e. there exists $g\in \kappa^\kappa$ such that $g\not \leq^*_\kappa f$ for every $f\in \maxk(Y)$.

For each $\xi< \kappa$, let $\mathcal{V}_{e(\xi)} = \{U^{(e(\xi))}_\eta : \eta< g(\xi)\}$. Then $|\mathcal{V}_{e(\xi)}| < \kappa$. Let $F\in [X]^{<\kappa}$. Define $h\in \maxk(Y)$ by $h(\xi) = \sup\{f_x(\xi) : x\in F\}$ for all $\xi< \kappa$, where $\sup \emptyset = 0$. Then there exists $\xi< \kappa$ such
that $g(\xi)> h(\xi)$. It follows that $f_x(\xi)< g(\xi)$ for every $x\in F$, and hence $F\subseteq \bigcup \mathcal{V}_{e(\xi)}$. Now enumerate $B = \{b_\xi : \xi< \kappa\}$. Recursively choose pairwise distinct $W_\xi \in \mathcal{U}_{b_\xi}$, and put $\mathcal{V}_{b_\xi} = \{W_\xi\}$. Then every $\mathcal{V}_\alpha$ has cardinality less than $\kappa$. Furthermore, since no $\mathcal{U}_\alpha$ has a subcover of cardinality less than $\kappa$, $\bigcup \mathcal{V}_\alpha \neq X$ for every $\alpha< \kappa$. The family $\{\bigcup \mathcal{V}_\alpha : \alpha< \kappa\}$ has cardinality $\kappa$. By the preceding argument, every member of $[X]^{<\kappa}$ is contained in some $\bigcup \mathcal{V}_\alpha$. Therefore, $\{\bigcup \mathcal{V}_\alpha : \alpha< \kappa\} \in \Omega_\kappa(X)$. Hence $X$ is $\kappa$-Scheepers.
\end{proof}

Using the preceding lemma, we prove only Theorem~\ref{thm6} concerning the $\kappa$-Scheepers property, as Theorems~\labelcref{thm5,thm7} can be proved similarly.

\begin{Th}
\label{thm5}
Assume that $\kappa$ is regular. Let $\mathcal{A} \subseteq [\kappa]^\kappa$ be an almost disjoint family. Then the following assertions are equivalent.
\begin{enumerate}[wide=0pt,label={\upshape(\arabic*)},leftmargin=*]
\item $\Psi_\kappa(\mathcal{A})$ is $\kappa$-Menger.
\item $\Psi_\kappa(\mathcal{A})$ is strongly star-$\kappa$-Menger.
\item $|\mathcal{A}|< \dk$.
\end{enumerate}
\end{Th}

\begin{Th}
\label{thm6}
Assume that $\kappa$ is regular. Let $\mathcal{A} \subseteq [\kappa]^\kappa$ be an almost disjoint family. Then the following assertions are equivalent.
\begin{enumerate}[wide=0pt,label={\upshape(\arabic*)},leftmargin=*]
\item $\Psi_\kappa(\mathcal{A})$ is $\kappa$-Scheepers.
\item $\Psi_\kappa(\mathcal{A})$ is strongly star-$\kappa$-Scheepers.
\item $|\mathcal{A}|< \dlk$.
\end{enumerate}
\end{Th}

\begin{proof}
We prove only the implication $(2)\Rightarrow (3)$, since all other implications follow immediately from Figure~\ref{fig1} and Lemma~\ref{lemma18}.

Assume that $|\mathcal{A}| \geq \dlk$. We show that $\Psi_\kappa(\mathcal{A})$ is not strongly star-$\kappa$-Scheepers. Choose a $\lessk$-dominating family $D\subseteq \kappa^\kappa$ of cardinality $\dlk$. Choose $\mathcal{B} \in [\mathcal{A}]^{\dlk}$ and enumerate $D$ as $D = \{f_B : B\in \mathcal{B}\}$. For $\alpha, \xi< \kappa$, let $\mathcal{B}_{\alpha,\xi} = \{B\in \mathcal{B} : f_B(\alpha) = \xi\}$. Define $\mathcal{P}_{\alpha,0} = (\mathcal{A} \setminus \mathcal{B}) \cup \mathcal{B}_{\alpha,0}$, and for $0< \xi< \kappa$, put $\mathcal{P}_{\alpha,\xi} = \mathcal{B}_{\alpha,\xi}$. Thus the nonempty members of $\{\mathcal{P}_{\alpha,\xi} : \xi< \kappa\}$ form a partition of $\mathcal{A}$.

For every $\alpha,\xi< \kappa$ such that $\mathcal{P}_{\alpha,\xi} \neq \emptyset$, define \[U_\xi^{(\alpha)} = \mathcal{P}_{\alpha,\xi} \cup \bigcup_{A\in \mathcal{P}_{\alpha,\xi}} (A\setminus \xi).\] Now put $\mathcal{U}_\alpha = \{U_\xi^{(\alpha)} : \xi< \kappa, \mathcal{P}_{\alpha,\xi} \neq \emptyset\} \cup \{\{\eta\} : \eta< \kappa\}$. Then $\mathcal{U}_\alpha \in \mathcal{O}_\kappa(\Psi_\kappa(\mathcal{A}))$ for every $\alpha< \kappa$. Let $F_\alpha \in [\Psi_\kappa(\mathcal{A})]^{<\kappa}$, $\alpha< \kappa$. Define $g\in \kappa^\kappa$ by \[g(\alpha) = \sup (\{\eta + 1 : \eta\in F_\alpha \cap \kappa\} \cup \{f_C(\alpha) + 1 : C\in F_\alpha \cap \mathcal{B}\} \cup \{1\}).\] Since $D$ is $\lessk$-dominating, we can pick a $\mathcal{C} \in [\mathcal{B}]^{<\kappa}$ such that for every $\alpha< \kappa$, there exists
$B_\alpha \in \mathcal{C}$ satisfying $g(\alpha) < f_{B_\alpha}(\alpha)$. Fix $\alpha< \kappa$ and put $\xi_\alpha = f_{B_\alpha}(\alpha)$. Since $g(\alpha) \geq 1$, we have $\xi_\alpha > 0$. We claim that $F_\alpha \cap U_{\xi_\alpha}^{(\alpha)} = \emptyset$. Let first $\eta\in F_\alpha \cap \kappa$. Then $\eta< g(\alpha) \leq \xi_\alpha$. Hence $\eta \notin \bigcup_{C\in \mathcal{P}_{\alpha,\xi_\alpha}} (C\setminus \xi_\alpha)$. Thus no isolated point of $F_\alpha$ belongs to $U_{\xi_\alpha}^{(\alpha)}$. Next let $C\in F_\alpha \cap \mathcal{B}$. By the definition of $g$, $f_C(\alpha) < g(\alpha) \leq \xi_\alpha = f_{B_\alpha}(\alpha)$, so $C\notin \mathcal{B}_{\alpha,\xi_\alpha} = \mathcal P_{\alpha,\xi_\alpha}$. Finally, if $C\in F_\alpha \cap (\mathcal{A} \setminus \mathcal{B})$, then $C\in \mathcal{P}_{\alpha,0}$, whereas $\xi_\alpha> 0$, and hence $C\notin \mathcal{P}_{\alpha,\xi_\alpha}$. Therefore, $F_\alpha \cap U_{\xi_\alpha}^{(\alpha)} = \emptyset$.

Clearly, $U_{\xi_\alpha}^{(\alpha)}$ is the only member of $\mathcal{U}_\alpha$ containing the non-isolated point $B_\alpha$. Consequently, $B_\alpha\notin St(F_\alpha, \mathcal{U}_\alpha)$. Since $\alpha< \kappa$ was arbitrary, $\mathcal{C}\nsubseteq \bigcup_{\alpha< \kappa} St(F_\alpha, \mathcal{U}_\alpha)$ for every $\alpha< \kappa$. Thus the sequence $(\mathcal{U}_\alpha : \alpha< \kappa)$ witnesses that $\Psi_\kappa(\mathcal{A})$ is not strongly star-$\kappa$-Scheepers.
\end{proof}

\begin{Th}
\label{thm7}
Assume that $\kappa$ is regular. Let $\mathcal{A} \subseteq [\kappa]^\kappa$ be an almost disjoint family. Then the following assertions are equivalent.
\begin{enumerate}[wide=0pt,label={\upshape(\arabic*)},leftmargin=*]
\item $\Psi_\kappa(\mathcal{A})$ is $\kappa$-Hurewicz.
\item $\Psi_\kappa(\mathcal{A})$ is strongly star-$\kappa$-Hurewicz.
\item $|\mathcal{A}|< \bk$.
\end{enumerate}
\end{Th}

We have already established that for meta-$\lessk$-compact spaces as well as for $\kappa$-paracompact spaces, the $\kappa$-Menger, $\kappa$-Scheepers, and $\kappa$-Hurewicz properties are equivalent to their corresponding strongly star variants (see Section~\ref{sie}). To illustrate that the preceding equivalences for generalized $\Psi$-spaces extend beyond these classes, we present an elementary observation (Observation~\ref{obs1}) followed by an example (Example~\ref{ex17}).

\begin{Obs}[{Folklore}]
\label{obs1}
Assume that $\kappa$ is regular. Let $\mathcal{A} \subseteq [\kappa]^\kappa$ be an almost disjoint family. Then the following assertions are equivalent.
\begin{enumerate}[wide=0pt,label={\upshape(\arabic*)},leftmargin=*]
\item $|\mathcal{A}| \leq \kappa$.
\item $\Psi_\kappa(\mathcal{A})$ is paracompact.
\item $\Psi_\kappa(\mathcal{A})$ is metacompact.
\item $\Psi_\kappa(\mathcal{A})$ is meta-$\lessk$-compact.
\end{enumerate}
\end{Obs}

The construction of the following example is based on the approach of Hru\v{s}\'{a}k et al. \cite{HMS} in the classical countable setting.

\begin{Ex}
\label{ex17}
Assume that $\kappa$ is regular and $\kappa^+ < \bk$. Then there exists an almost disjoint family $\mathcal{A} \subseteq [\kappa]^\kappa$ such that $|\mathcal{A}| < \bk$ and $\Psi_\kappa(\mathcal{A})$ is not $\kappa$-paracompact.
\end{Ex}

\begin{proof}
We first construct an almost disjoint family $\mathcal{A} = \{A_\alpha : \alpha< \kappa^+\}$ such that the size of $\{\beta< \alpha : A_\alpha \cap A_\beta \subseteq \delta\}$ is less than $\kappa$ for every $\kappa \leq \alpha< \kappa^+$ and $\delta< \kappa$. Choose pairwise disjoint sets $A_\alpha\in [\kappa]^\kappa$ for $\alpha< \kappa$. Suppose that $\kappa \leq \alpha< \kappa^+$ and that $A_\beta$ has been constructed for every $\beta< \alpha$. Enumerate the preceding family as \[\{A_\beta : \beta< \alpha\} = \{B_\xi : \xi< \kappa\}.\] For each $\xi< \kappa$, choose \[x_\xi\in B_\xi \setminus \left((\xi+1) \cup \bigcup_{\eta< \xi} B_\eta\right).\] This is possible because $\kappa$ is regular, the sets $B_\xi$ are pairwise almost disjoint, and $|B_\xi| = \kappa$. Put $A_\alpha = \{x_\xi : \xi< \kappa\}$. Clearly, $|A_\alpha| = \kappa$. Moreover, for each $\xi< \kappa$, $A_\alpha \cap B_\xi \subseteq \{x_\eta : \eta\leq \xi\}$, and so $|A_\alpha \cap B_\xi|< \kappa$. Thus the family remains almost disjoint. Furthermore, $x_\xi> \xi$ and $x_\xi\in A_\alpha \cap B_\xi$. Consequently, for every $\delta< \kappa$, all but fewer than $\kappa$ preceding members $B_\xi$ satisfy $A_\alpha \cap B_\xi\not \subseteq \delta$. This completes the construction.

Now partition $\kappa^+$ into countably many pairwise disjoint sets $I_n$, $n\in \omega$, each of cardinality $\kappa^+$. For every $n\in \omega$, put $U_n = \kappa \cup \{A_\alpha : \alpha\in I_n\}$. Then $\mathcal{U} = \{U_n : n\in \omega\}$ is a countable open cover of $\Psi_\kappa(\mathcal{A})$. Suppose that $\mathcal{U}$ admits a locally finite open refinement $\mathcal{V}$. By grouping the members of $\mathcal{V}$ according to the corresponding members of $\mathcal{U}$, we obtain a locally finite open cover $\{W_n : n\in \omega\}$ such that $W_n\subseteq U_n$ for every $n\in \omega$.

For each $\alpha< \kappa^+$, let $n_\alpha$ be the unique $n\in \omega$ such that $\alpha\in I_n$. Since $A_\alpha$ belongs only to $U_{n_\alpha}$, we have $A_\alpha\in W_{n_\alpha}$. Choose $\delta_\alpha < \kappa$ such that $N(A_\alpha, \delta_\alpha) \subseteq W_{n_\alpha}$. For each $n\in \omega$, since $|I_n| = \kappa^+$, there exist $\delta_n< \kappa$ and $H_n\subseteq I_n$ with $|H_n| = \kappa$ such that $\delta_\beta = \delta_n$ for every $\beta\in H_n$. Since $|\bigcup_{n\in \omega} H_n| = \kappa< \kappa^+$, we may choose $\alpha< \kappa^+$ such that $\bigcup_{n\in \omega} H_n \subseteq \alpha$.

By the local finiteness of $\{W_n : n\in \omega\}$, there exists $\gamma< \kappa$ such that $N(A_\alpha, \gamma)$ meets only finitely many $W_n$. Fix $n\in \omega$ and put $\delta = \max\{\gamma, \delta_n\}$. By the construction of $\mathcal{A}$, there exists $\beta\in H_n$ such that $A_\alpha \cap A_\beta \not \subseteq \delta$. Choose $\xi\in (A_\alpha \cap A_\beta) \setminus \delta$. Then \[\xi\in N(A_\alpha, \gamma) \cap N(A_\beta, \delta_n) \subseteq N(A_\alpha, \gamma) \cap W_n.\] Since $n\in \omega$ was arbitrary, $N(A_\alpha, \gamma)$ meets every $W_n$, contradicting local finiteness. Therefore, $\Psi_\kappa(\mathcal{A})$ is not $\kappa$-paracompact.

Finally, $|\mathcal{A}| = \kappa^+< \bk$, as required.
\end{proof}

We now turn to the $\kappa$-Rothberger property. The proof of the following lemma is analogous to that of Lemma~\ref{lemma18} and is therefore omitted.

\begin{Lemma}
\label{lemma21}
Assume that $\kappa^{<\kappa} = \kappa$. If $|X|< \ink$, then $X$ is $\kappa$-Rothberger.
\end{Lemma}

In contrast to the classical countable case, the preceding lemma yields the following result.

\begin{Th}
\label{thm8}
Assume that $\kappa^{<\kappa} = \kappa$. Let $\mathcal{A}\subseteq [\kappa]^\kappa$ be an almost disjoint family. Then $\Psi_\kappa(\mathcal{A})$ is $\kappa$-Rothberger if and only if $|\mathcal{A}|< \ink$.
\end{Th}

\begin{proof}
If $|\mathcal{A}| < \ink$, then by Lemma~\ref{lemma21}, $\Psi_\kappa(\mathcal{A})$ is $\kappa$-Rothberger. For the other implication, suppose that $|\mathcal{A}| \geq \ink$. Choose an eventually different family $E\subseteq \kappa^\kappa$ of cardinality $\ink$. For each $f\in E$ and $\beta, \xi< \kappa$, define $f_{\beta,\xi} \in \kappa^\kappa$ by
\[f_{\beta,\xi}(\alpha)= \begin{cases}
\xi, & \text{if } \alpha< \beta,\\
f(\alpha), & \text{if } \alpha \geq \beta.
\end{cases}\]
Put $E^\prime = \{f_{\beta,\xi} : f\in E \text{ and } \beta, \xi< \kappa\}$. Since $\kappa< \ink$, we have $|E^\prime| = \ink$. Moreover, for every $g\in \kappa^\kappa$, there exists $h\in E^\prime$ such that $h(\alpha) \neq g(\alpha)$ for every $\alpha< \kappa$. Indeed, choose $f\in E$ and $\beta< \kappa$ such that $f(\alpha) \neq g(\alpha)$ for every $\alpha \geq \beta$, and choose $\xi\in \kappa \setminus g[\beta]$. Then $h = f_{\beta,\xi}$ has the desired property.

Choose a surjection $A\mapsto f_A$ from $\mathcal{A}$ onto $E^\prime$. For each $\alpha, \xi< \kappa$, put $U_\xi^{(\alpha)} = \kappa \cup \{A\in \mathcal{A} : f_A(\alpha) = \xi\}$. For each $\alpha< \kappa$, let $\mathcal{U}_\alpha = \{U_\xi^{(\alpha)} : \xi< \kappa\} \cup \{\{\xi\} : \xi< \kappa\}$. Clearly, $\mathcal{U}_\alpha \in \mathcal{O}_\kappa(\Psi_\kappa(\mathcal{A}))$. Let $U_\alpha \in \mathcal{U}_\alpha$, $\alpha< \kappa$, be arbitrary. Define $g\in \kappa^\kappa$ by choosing $g(\alpha)$ such that $U_\alpha = U_{g(\alpha)}^{(\alpha)}$ whenever possible, and putting $g(\alpha) = 0$ otherwise. By the choice of $E^\prime$, there exists $A\in \mathcal{A}$ such that $f_A(\alpha) \neq g(\alpha)$ for every $\alpha< \kappa$. Consequently, $A\notin \bigcup_{\alpha< \kappa} U_\alpha$. Thus $\Psi_\kappa(\mathcal{A})$ is not $\kappa$-Rothberger.
\end{proof}

\begin{Cor}
\label{cor14}
Assume that $\kappa^{<\kappa} = \kappa$. Let $\mathcal{A} \subseteq [\kappa]^\kappa$ be an almost disjoint family. If $|\mathcal{A}|< \ink$, then $\Psi_\kappa(\mathcal{A})$ is strongly star-$\kappa$-Rothberger.
\end{Cor}

Unlike the corresponding results for the strongly star-$\kappa$-Menger, strongly star-$\kappa$-Scheepers, and strongly star-$\kappa$-Hurewicz properties, the preceding corollary provides only a sufficient cardinality condition for the strongly star-$\kappa$-Rothberger property of generalized $\Psi$-spaces. The following example shows that this cardinality bound is sharp.

\begin{Ex}
\label{ex15}
Assume that $\kappa^{<\kappa} = \kappa$. Then there exists an almost disjoint family $\mathcal{A} \subseteq [\kappa]^\kappa$ of cardinality $\ink$ such that $\Psi_\kappa(\mathcal{A})$ is not star-$\kappa$-Rothberger. Consequently, $\Psi_\kappa(\mathcal{A})$ is not strongly star-$\kappa$-Rothberger.
\end{Ex}

\begin{proof}
There exists a subspace $X\subseteq \kappa^\kappa$ such that $|X| = \ink$ and $X$ is not $\kappa$-Rothberger (see \cite{AC}). Then there exists a sequence $(\mathcal{U}_\alpha : \alpha< \kappa)\in \mathcal{O}_\kappa(X)^\kappa$, where $\mathcal{U}_\alpha = \{U_\xi^{(\alpha)} : \xi< \kappa\}$, such that for every $g\in \kappa^\kappa$, $X\neq \bigcup_{\alpha< \kappa} U_{g(\alpha)}^{(\alpha)}$. For each $x\in X$, define $s_x\in \kappa^\kappa$ by \[s_x(\alpha) = \min \{\xi< \kappa : x\in U_\xi^{(\alpha)}\}.\] Fix a bijection $\pi: \kappa \times \kappa\to \kappa$, and define $h_x\in \kappa^\kappa$ by $h_x(\alpha) = \pi(s_x(\alpha), x(\alpha))$ for every $\alpha< \kappa$. Then the family $B = \{h_x : x\in X\}$ has cardinality $\ink$. We claim that for each $g\in \kappa^\kappa$, there exists $h\in B$ such that $h(\alpha) \neq g(\alpha)$ for every $\alpha< \kappa$. To see this, let $g\in \kappa^\kappa$. For each $\alpha< \kappa$, write $\pi^{-1}(g(\alpha)) = (\xi_\alpha, \eta_\alpha)$. Then there exists $x\in X$ such that $x\notin U_{\xi_\alpha}^{(\alpha)}$ for every $\alpha< \kappa$. Hence $s_x(\alpha) \neq \xi_\alpha$ for every $\alpha< \kappa$, and therefore \[h_x(\alpha) = \pi(s_x(\alpha), x(\alpha)) \neq \pi(\xi_\alpha, \eta_\alpha) = g(\alpha)\] for every $\alpha< \kappa$. This proves the claim.

Partition $\kappa$ into two disjoint sets $I$ and $J$ with $|I| = |J| = \kappa$, and fix a bijection $\varphi: \kappa^{<\kappa} \to I$. For each $h\in B$, put $A_h = \{\varphi(h\restriction (\beta+1)) : \beta< \kappa\}$. Clearly, $|A_h| = \kappa$. Moreover, if $h, h^\prime\in B$ are distinct and $\delta< \kappa$ is the least ordinal such that $h(\delta) \neq h^\prime(\delta)$, then $|A_h\cap A_{h^\prime}| \leq |\delta| < \kappa$. Thus $\mathcal{A} = \{A_h : h\in B\}$ is an almost disjoint family with $|\mathcal{A}| = \ink$.

We now show that $\Psi_\kappa(\mathcal{A})$ is not star-$\kappa$-Rothberger. For each $h\in B$ and each $\alpha< \kappa$, by the regularity of $\kappa$, choose
$\gamma_{h, \alpha} < \kappa$ such that $\varphi(h\restriction (\beta+1)) < \gamma_{h, \alpha}$ whenever $\beta\leq \alpha$. For every $\alpha, \xi< \kappa$, define \[V_\xi^{(\alpha)} = \bigcup_{\substack{h\in B\\ h(\alpha) = \xi}} (\{A_h\} \cup (A_h \setminus \gamma_{h, \alpha})).\] Each $V_\xi^{(\alpha)}$ is open in $\Psi_\kappa(\mathcal{A})$. We claim that for each fixed $\alpha< \kappa$, the nonempty sets $V_\xi^{(\alpha)}$, $\xi< \kappa$, are pairwise disjoint. Indeed, suppose that $z\in V_\xi^{(\alpha)} \cap V_\eta^{(\alpha)}$. If $z = A_h\in \mathcal{A}$, then $\xi = h(\alpha) = \eta$. Otherwise, $z\in \kappa$. Then there exist $h, h^\prime\in B$ and $\beta, \delta> \alpha$ such that $z = \varphi(h\restriction (\beta+1)) = \varphi(h^\prime \restriction (\delta+1))$, where $h(\alpha) = \xi$ and $h^\prime(\alpha) = \eta$. It follows that $h\restriction (\beta+1) = h^\prime \restriction(\delta+1)$. In particular, $h(\alpha) = h^\prime(\alpha)$ and so $\xi = \eta$. This proves the claim.

For each $\alpha< \kappa$, let \[\mathcal{V}_\alpha = \{V_\xi^{(\alpha)} : \xi< \kappa, V_\xi^{(\alpha)} \neq \emptyset\} \cup \{\{\zeta\} : \zeta< \kappa\}.\] Then $\mathcal{V}_\alpha$ is an open cover of $\Psi_\kappa(\mathcal{A})$ of cardinality $\kappa$. Moreover, it has no subcover of cardinality less than $\kappa$, since every point of $J$ belongs only to its singleton member of $\mathcal{V}_\alpha$. Then the sequence $(\mathcal{V}_\alpha : \alpha< \kappa) \in \mathcal{O}_\kappa(\Psi_\kappa(\mathcal{A}))^\kappa$ witnesses that $\Psi_\kappa(\mathcal{A})$ is not star-$\kappa$-Rothberger. Indeed, for each $\alpha< \kappa$, select $W_\alpha \in \mathcal{V}_\alpha$. Define $g\in \kappa^\kappa$ by
\[g(\alpha) = \begin{cases}
\xi, & \mbox{if } W_\alpha \cap V_\xi^{(\alpha)} \neq \emptyset \text{ for some } \xi< \kappa\\
g(\alpha) = 0, & \mbox{otherwise}.
\end{cases}\]
Now let $h\in B$ and $\alpha< \kappa$. Since $V_{h(\alpha)}^{(\alpha)}$ is the unique member of $\mathcal{V}_\alpha$ containing $A_h$, we have $A_h\in St(W_\alpha, \mathcal{V}_\alpha)$ only if $W_\alpha \cap V_{h(\alpha)}^{(\alpha)} \neq \emptyset$. Consequently, $A_h\in St(W_\alpha, \mathcal{V}_\alpha)$ gives $h(\alpha) = g(\alpha)$. Which is absurd. Thus $A_h\notin St(W_\alpha, \mathcal{V}_\alpha)$ for every $\alpha< \kappa$. Therefore, $\Psi_\kappa(\mathcal{A})$ is not star-$\kappa$-Rothberger.
\end{proof}

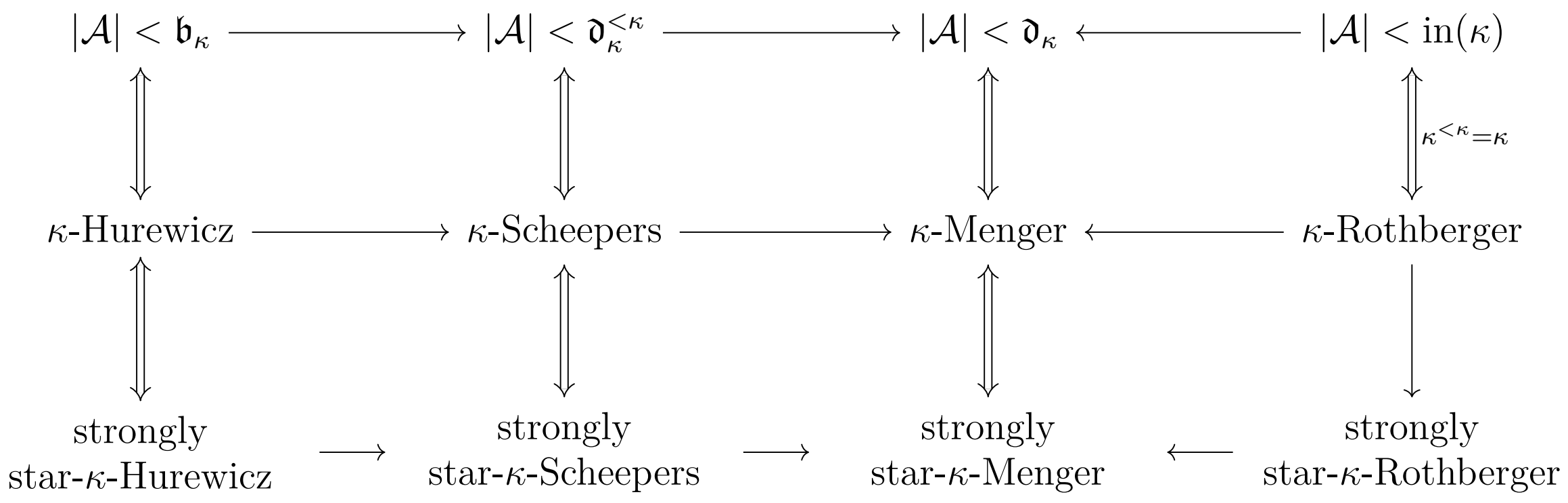
\begin{figure}[h]
\begin{adjustbox}{max width=\textwidth,max height=\textheight,keepaspectratio,center}
\begin{tikzcd}[column sep=4ex, row sep=8ex, arrows={crossing over}]
\lvert\mathcal{A}\rvert < \bk \arrow[r] & \lvert\mathcal{A}\rvert < \dlk \arrow[r] & \lvert\mathcal{A}\rvert < \dk & \lvert\mathcal{A}\rvert < \ink \arrow[l] \arrow[d, Leftrightarrow, "\kappa^{<\kappa} = \kappa"]\\
\text{$\kappa$-Hurewicz} \arrow[r] \arrow[u, Leftrightarrow] & \text{$\kappa$-Scheepers} \arrow[r] \arrow[u, Leftrightarrow] & \text{$\kappa$-Menger} \arrow[u, Leftrightarrow] & \text{$\kappa$-Rothberger} \arrow[l] \arrow[d]\\
\text{\parbox{3.5cm}{\centering strongly\\ star-$\kappa$-Hurewicz}} \arrow[r] \arrow[u, Leftrightarrow] & \text{\parbox{3.5cm}{\centering strongly\\ star-$\kappa$-Scheepers}} \arrow[r] \arrow[u, Leftrightarrow] & \text{\parbox{3.5cm}{\centering strongly\\ star-$\kappa$-Menger}} \arrow[u, Leftrightarrow] & \text{\parbox{3.5cm}{\centering strongly\\ star-$\kappa$-Rothberger}} \arrow[l]
\end{tikzcd}
\end{adjustbox}
\caption{Cardinality characterizations and implications for generalized $\Psi$-spaces.}
\label{fig3}
\end{figure}

The cardinality characterizations and implications obtained for generalized $\Psi$-spaces are summarized in Figure~\ref{fig3}, where $\kappa$ is assumed to be regular.

\begin{Rem}
\label{rem1}
The reverse implications along the rows of Figure~\ref{fig3} do not hold in general, as can be seen by imposing suitable cardinality restrictions on the underlying almost disjoint families under appropriate inequalities between the relevant cardinal invariants. However, the remaining reverse implication along the fourth column remains open.
\end{Rem}

\subsection{Failure of star-$\kappa$ selection properties}
The following observation follows from elementary cardinal arithmetic (cf. \cite[Lemma~8.1]{CFHT}). We include a sketch of the proof for the convenience of the reader.

\begin{Lemma}
\label{lemma25}
Let $\lambda$ be an infinite cardinal such that $\lambda^\kappa = \lambda$. Then $\cof([\lambda]^\kappa) = \lambda$.
\end{Lemma}

\begin{proof}
Clearly, $\cof([\lambda]^\kappa) \leq \lambda$. Conversely, let $\mathcal{D} \subseteq [\lambda]^\kappa$ be a cofinal family. Then $\bigcup \mathcal{D} = \lambda$. Suppose that $|\mathcal{D}|< \lambda$. Since every member of $\mathcal{D}$ has cardinality $\kappa$ and $\kappa< \lambda$, we have $\lambda = |\bigcup \mathcal{D}| \leq| \mathcal{D}| \cdot \kappa< \lambda$, a contradiction. Therefore, $|\mathcal{D}| \geq \lambda$, and hence $\cof([\lambda]^\kappa) \geq \lambda$. Combining the two inequalities, we conclude that $\cof([\lambda]^\kappa) = \lambda$.
\end{proof}

We now observe that the natural uncountable analogues of \cite[Theorems~1.9 and 2.4]{tsaban15} fail in general.

\begin{Ex}
\label{ex16}
Assume that $\kappa^{<\kappa} = \kappa$. Then there exists an almost disjoint family $\mathcal{A} \subseteq [\kappa]^\kappa$ of cardinality $\lambda$ such that $\lambda \geq \dk$ and $\cof([\lambda]^\kappa) = \lambda$, but $\Psi_\kappa(\mathcal{A})$ is star-$\kappa$-Hurewicz.
\end{Ex}

\begin{proof}
Since $\kappa^{<\kappa} = \kappa$, it follows from \cite[Observation~2.2]{ms20} that there exists an almost disjoint family $\mathcal{A} \subseteq [\kappa]^\kappa$ of cardinality $2^\kappa$. Extend this family to a MAD family $\mathcal{A}$. Clearly, $|\mathcal{A}| = 2^\kappa$.

We show that for every $\mathcal{U} \in \mathcal{O}_\kappa(\Psi_\kappa(\mathcal{A}))$, there exists $\mathcal{V} \in [\mathcal{U}]^{<\kappa}$ such that $\Psi_\kappa(\mathcal{A}) = St(\bigcup \mathcal{V}, \mathcal{U})$. Let $\mathcal{U} = \{U_\alpha : \alpha< \kappa\}\in \mathcal{O}_\kappa(\Psi_\kappa(\mathcal{A}))$. We first claim that there exists
$\mathcal{V} \in [\mathcal{U}]^{<\kappa}$ satisfying $\kappa \subseteq \bigcup \mathcal{V}$. Suppose, to the contrary, that no subfamily of $\mathcal{U}$ of cardinality less than $\kappa$ covers $\kappa$. Then for every $\alpha< \kappa$, $|\kappa \setminus \bigcup_{\beta \leq \alpha} U_\beta| = \kappa$. By transfinite recursion, choose distinct points $x_\alpha\in \kappa$, $\alpha< \kappa$, such that $x_\alpha \notin \bigcup_{\beta \leq \alpha} U_\beta$. Put $D = \{x_\alpha : \alpha< \kappa\}$. Then $|D| = \kappa$. By the maximality of $\mathcal{A}$, there exists $A\in \mathcal{A}$ such that $|A\cap D| = \kappa$. Choose $\beta< \kappa$ with $A\in U_\beta$. Since $U_\beta$ is an open neighbourhood of $A$ in $\Psi_\kappa(\mathcal{A})$, we have $|A\setminus U_\beta| < \kappa$. Consequently, $|D\cap U_\beta| = \kappa$. On the other hand, by the construction of $D$, $D\cap U_\beta \subseteq \{x_\alpha : \alpha< \beta\}$, and hence $|D\cap U_\beta| < \kappa$, a contradiction. Therefore, there exists $\mathcal{V}\in [\mathcal{U}]^{<\kappa}$ such that $\kappa \subseteq \bigcup \mathcal{V}$. Thus $\Psi_\kappa(\mathcal{A}) = St(\bigcup \mathcal{V}, \mathcal{U})$. It follows immediately that $\Psi_\kappa(\mathcal{A})$ is star-$\kappa$-Hurewicz.

Finally, let $\lambda = 2^\kappa$. Then $\lambda^\kappa = (2^\kappa)^\kappa = 2^\kappa = \lambda$ and $\lambda \geq \dk$. By Lemma~\ref{lemma25}, $\cof([\lambda]^\kappa) = \lambda$. This completes the proof.
\end{proof}

However, we now present constructions of generalized $\Psi$-spaces that fail to be star-$\kappa$-Menger, star-$\kappa$-Scheepers, or star-$\kappa$-Hurewicz; see Theorems~\labelcref{thm9,thm10,thm11}, respectively.  For functions $s$ and $t$ whose domains are ordinals, we write $s\sqsubseteq t$ if $s = t\restriction \dom(s)$, i.e. if $s$ is an initial segment of $t$.

\begin{Th}
\label{thm9}
Assume that $\kappa^{<\kappa} = \kappa$. Let $e: \kappa^{<\kappa} \to \kappa$ be a bijection and $D\subseteq \kappa^\kappa$ be pointwise dominating. For each $f\in D$, put $A_f = \{e(f\restriction \xi) : \xi< \kappa\}$. Then $\mathcal{A} = \{A_f : f\in D\}$ is an almost disjoint family and $\Psi_\kappa(\mathcal{A})$ is not star-$\kappa$-Menger.
\end{Th}

\begin{proof}
Clearly, $|A_f| = \kappa$ for every $f\in D$. If $f,g\in D$ are distinct, then there exists $\xi< \kappa$ such that $f(\xi) \neq g(\xi)$. Consequently, $A_f\cap A_g \subseteq \{e(f\restriction \eta) : \eta \leq \xi\}$, and hence $|A_f\cap A_g|< \kappa$. Thus $\mathcal{A}$ is an almost disjoint
family.

For each $s\in\kappa^{<\kappa}$, define \[W_s = \{A_f : f\in D, s\sqsubseteq f\} \cup \{e(t) : t\in \kappa^{<\kappa}, s\sqsubseteq t\}.\] Then $W_s$ is open in $\Psi_\kappa(\mathcal{A})$. For each $\alpha< \kappa$, let \[\mathcal{U}_\alpha = \{W_s : s\in \kappa^{\alpha+1}\} \cup \{\{e(t)\} : t\in \kappa^{<\kappa}\}.\] Then $\mathcal{U}_\alpha \in \mathcal{O}_\kappa(\Psi_\kappa(\mathcal{A}))$. Let $\mathcal{V}_\alpha \in [\mathcal{U}_\alpha]^{<\kappa}$, $\alpha< \kappa$, be arbitrary. For every $\alpha< \kappa$, define \[E_\alpha = \{s\in \kappa^{\alpha+1} : W_s\in \mathcal{V}_\alpha\} \cup \{t\restriction(\alpha+1) : t\in \kappa^{<\kappa}, \dom(t) \geq \alpha+1, \{e(t)\} \in \mathcal{V}_\alpha\}.\] Clearly, $|E_\alpha| < \kappa$. Define $g\in \kappa^\kappa$ by \[g(\alpha) = \sup (\{s(\alpha)+1 : s\in E_\alpha\} \cup \{0\})\] for every $\alpha< \kappa$. Since $D$ is pointwise dominating, there exists $f\in D$ such that $g(\alpha) \leq f(\alpha)$ for every $\alpha< \kappa$. It follows that $f\restriction (\alpha+1) \notin E_\alpha$ for every $\alpha< \kappa$. Fix $\alpha< \kappa$. Notice that $W_{f\restriction (\alpha+1)}$ is the unique member of $\mathcal{U}_\alpha$ containing $A_f$. Indeed, for distinct $s,t\in \kappa^{\alpha+1}$, the sets $W_s$ and $W_t$ are disjoint.

We claim that $W_{f\restriction (\alpha+1)} \cap \bigcup \mathcal{V}_\alpha = \emptyset$. To see this, first suppose that $W_t\in \mathcal{V}_\alpha$ for some $t\in \kappa^{\alpha+1}$. Then $t\in E_\alpha$ and so $t\neq f\restriction (\alpha+1)$. Consequently, $W_{f\restriction (\alpha+1)} \cap W_t = \emptyset$. Next, suppose that $\{e(t)\} \in \mathcal{V}_\alpha$ for some $t\in \kappa^{<\kappa}$. If $\dom(t)< \alpha+1$, then $e(t)\notin W_{f\restriction (\alpha+1)}$. Otherwise, $t\restriction (\alpha+1) \in E_\alpha$, and hence $t\restriction (\alpha+1) \neq f\restriction (\alpha+1)$. Again, $e(t)\notin W_{f\restriction (\alpha+1)}$. Therefore, $W_{f\restriction (\alpha+1)} \cap \bigcup \mathcal{V}_\alpha = \emptyset$.

Since $W_{f\restriction (\alpha+1)}$ is the unique member of $\mathcal{U}_\alpha$ containing $A_f$, we obtain $A_f\notin St(\bigcup \mathcal{V}_\alpha, \mathcal{U}_\alpha)$. As $\alpha< \kappa$ was arbitrary, $A_f\notin \bigcup_{\alpha< \kappa} St(\bigcup \mathcal{V}_\alpha, \mathcal{U}_\alpha)$. Thus the sequence $(\mathcal{U}_\alpha : \alpha< \kappa)$ witnesses that $\Psi_\kappa(\mathcal{A})$ is not star-$\kappa$-Menger.
\end{proof}

\begin{Cor}
\label{cor3}
Assume that $\kappa^{<\kappa} = \kappa$. Then for every cardinal $\lambda$ with $\dk\leq \lambda \leq 2^\kappa$, there exists an almost disjoint family $\mathcal{A} \subseteq [\kappa]^\kappa$ of cardinality $\lambda$ such that $\Psi_\kappa(\mathcal{A})$ is not star-$\kappa$-Menger.
\end{Cor}

\begin{proof}
Let $D\subseteq \kappa^\kappa$ be a dominating family of cardinality
$\dk$. For each $d\in D$, $\xi< \kappa$, and $h\in \kappa^\xi$, define $d_{\xi,h} \in \kappa^\kappa$ by
\[d_{\xi,h}(\alpha) = \begin{cases}
h(\alpha), & \text{if } \alpha< \xi,\\
d(\alpha), & \text{if } \alpha \geq \xi.
\end{cases}\]
Put $A = \{d_{\xi,h} : d\in D, \xi< \kappa, h\in \kappa^\xi\}$. Since $\kappa^{<\kappa} = \kappa$, we have $|A| = \dk$. We claim that $A$ is pointwise dominating. Indeed, let $g\in \kappa^\kappa$. Choose $d\in D$ such that $g\leq^*_\kappa d$. Then there exists $\xi< \kappa$ such that $g(\alpha) \leq d(\alpha)$ for all $\alpha \geq \xi$. Taking $h = g\restriction \xi$, we obtain $g(\alpha) \leq d_{\xi,h}(\alpha)$ for every $\alpha< \kappa$. Thus $A$ is pointwise dominating.

Now let $\lambda$ satisfy $\dk\leq \lambda \leq 2^\kappa$. Since $|\kappa^\kappa| = 2^\kappa$, choose $B \subseteq \kappa^\kappa$ such that $A \subseteq B$ and $|B| = \lambda$. Fix a bijection $e: \kappa^{<\kappa} \to \kappa$, and for every $f\in B$, define $A_f = \{e(f\restriction \xi) : \xi< \kappa\}$. Let $\mathcal{A} = \{A_f : f\in B\}$. Then $|\mathcal{A}| = \lambda$. By Theorem~\ref{thm9}, $\mathcal{A}$ is an almost disjoint family and $\Psi_\kappa(\mathcal{A})$ is not star-$\kappa$-Menger.
\end{proof}

The proofs of Theorems~\ref{thm10} and \ref{thm11} are analogous to that of Theorem~\ref{thm9}, with suitable modifications, and are therefore removed.

\begin{Th}
\label{thm10}
Assume that $\kappa^{<\kappa} = \kappa$. Let $e: \kappa^{<\kappa} \to \kappa$ be a bijection and $D\subseteq \kappa^\kappa$ be a $\lessk$-dominating family. For each $f\in D$, put $A_f = \{e(f\restriction \xi) : \xi< \kappa\}$. Then $\mathcal{A} = \{A_f : f\in D\}$ is an almost disjoint family and $\Psi_\kappa(\mathcal{A})$ is not star-$\kappa$-Scheepers.
\end{Th}

\begin{Cor}
\label{cor4}
Assume that $\kappa^{<\kappa} = \kappa$. Then for every cardinal $\lambda$ with $\dlk \leq \lambda \leq 2^\kappa$, there exists an almost disjoint family $\mathcal{A} \subseteq [\kappa]^\kappa$ of cardinality $\lambda$ such that $\Psi_\kappa(\mathcal{A})$ is not star-$\kappa$-Scheepers.
\end{Cor}

\begin{Th}
\label{thm11}
Assume that $\kappa^{<\kappa} = \kappa$. Let $e: \kappa^{<\kappa} \to \kappa$ be a bijection and $B\subseteq \kappa^\kappa$ be an unbounded family. For each $f\in B$, put $A_f = \{e(f\restriction \xi) : \xi< \kappa\}$. Then $B = \{A_f : f\in B\}$ is an almost disjoint family and $\Psi_\kappa(\mathcal{A})$ is not star-$\kappa$-Hurewicz.
\end{Th}

\begin{Cor}
\label{cor5}
Assume that $\kappa^{<\kappa} = \kappa$. Then for every cardinal $\lambda$ with $\bk\leq \lambda \leq 2^\kappa$, there exists an almost disjoint family $\mathcal{A} \subseteq [\kappa]^\kappa$ of cardinality $\lambda$ such that $\Psi_\kappa(\mathcal{A})$ is not star-$\kappa$-Hurewicz.
\end{Cor}

In \cite[Problem~4.4]{tsaban15}, Tsaban posed the question of whether an analogue of \cite[Theorems~1.9 and 2.4]{tsaban15} exists for the star-Rothberger property (a space $X$ is said to be star-Rothberger \cite{kocinac99} if for every sequence $(\mathcal{U}_n : n\in \omega)$ of open covers of $X$, there exist $U_n\in \mathcal{U}_n$, $n\in \omega$, such that $X = \bigcup_{n\in \omega} St(U_n, \mathcal{U}_n)$). The corresponding uncountable analogue, as stated below, also remains open.

\begin{Prob}
\label{prob1}
Assume that $\kappa^{<\kappa} = \kappa$. Let $\mathcal{A} \subseteq [\kappa]^\kappa$ be an almost disjoint family of cardinality
$\lambda \geq \ink$. If $\cof([\lambda]^\kappa) = \lambda$, must $\Psi_\kappa(\mathcal{A})$ fail to be star-$\kappa$-Rothberger?
\end{Prob}

We now present a construction of generalized $\Psi$-space that fails to be star-$\kappa$-Rothberger. The proof follows an argument similar to that of Theorem~\ref{thm9}, with necessary modifications. For the convenience of the reader, we provide a detailed proof.

\begin{Th}
\label{thm12}
Assume that $\kappa^{<\kappa} = \kappa$. Let $e: \kappa^{<\kappa} \to \kappa$ be a bijection and $E\subseteq \kappa^\kappa$ satisfy the following condition: \[(\forall g\in \kappa^\kappa)\ (\exists f\in E)\ (\forall \alpha< \kappa)\ (f(\alpha) \neq g(\alpha)).\] For each $f\in E$, put $A_f = \{e(f\restriction \xi) : \xi< \kappa\}$. Then $\mathcal{A} = \{A_f : f\in E\}$ is an almost disjoint family and $\Psi_\kappa(\mathcal{A})$ is not star-$\kappa$-Rothberger.
\end{Th}

\begin{proof}
As in the proof of Theorem~\ref{thm9}, $\mathcal{A}$ is an almost disjoint family. For each $\alpha< \kappa$, \[\mathcal{U}_\alpha = \{W_s : s\in \kappa^{\alpha+1}\} \cup \{\{e(t)\} : t\in \kappa^{<\kappa}\} \in \mathcal{O}_\kappa(\Psi_\kappa(\mathcal{A})),\] where $W_s = \{A_f : f\in E, s\sqsubseteq f\} \cup \{e(t) : t\in \kappa^{<\kappa}, s\sqsubseteq t\}$ for every $s\in\kappa^{<\kappa}$.  Let $U_\alpha \in \mathcal{U}_\alpha$, $\alpha< \kappa$, be arbitrary. Define $g\in \kappa^\kappa$ by
\[g(\alpha) = \begin{cases}
s(\alpha), & \text{if } U_\alpha = W_s \text{ for some } s\in \kappa^{\alpha+1},\\
t(\alpha), & \text{if } U_\alpha = \{e(t)\} \text{ for some } t\in \kappa^{<\kappa} \text{ with } \dom(t)> \alpha,\\
0, & \text{otherwise}.
\end{cases}\]
By our assumption on $E$, there exists $f\in E$ such that $f(\alpha) \neq g(\alpha)$ for every $\alpha< \kappa$. We claim that $A_f\notin St(U_\alpha, \mathcal{U}_\alpha)$ for every $\alpha< \kappa$. Fix $\alpha< \kappa$. Observe that $W_{f\restriction (\alpha+1)} \cap U_\alpha = \emptyset$. Indeed, suppose first that $U_\alpha = W_s$ for some $s\in \kappa^{\alpha+1}$. By the definition of $g$, $f(\alpha) \neq s(\alpha)$. Thus $f\restriction (\alpha+1)\neq s$. Since $f\restriction (\alpha+1)$ and $s$ have the same length, they are incompatible. Therefore, $W_{f\restriction (\alpha+1)} \cap U_\alpha = W_{f\restriction (\alpha+1)} \cap W_s = \emptyset$. Next, suppose that $U_\alpha = \{e(t)\}$ for some $t\in \kappa^{<\kappa}$. If $\dom(t) \leq \alpha$, then $f\restriction (\alpha+1)\not \sqsubseteq t$, and hence $e(t)\notin W_{f\restriction (\alpha+1)}$. If $\dom(t) > \alpha$, then $f(\alpha) \neq g(\alpha) = t(\alpha)$, and again $f\restriction (\alpha+1)\not \sqsubseteq t$. Consequently, $W_{f\restriction (\alpha+1)} \cap U_\alpha = \emptyset$.

Hence, in either case, $W_{f\restriction (\alpha+1)} \cap U_\alpha = \emptyset$. Since $W_{f\restriction (\alpha+1)}$ is the unique member of $\mathcal{U}_\alpha$ containing $A_f$, $A_f\notin St(U_\alpha, \mathcal{U}_\alpha)$. Thus $A_f\notin \bigcup_{\alpha< \kappa} St(U_\alpha, \mathcal{U}_\alpha)$. Therefore, the sequence $(\mathcal{U}_\alpha : \alpha< \kappa)$ witnesses that $\Psi_\kappa(\mathcal{A})$ is not star-$\kappa$-Rothberger.
\end{proof}

\begin{Cor}
\label{cor6}
Assume that $\kappa^{<\kappa} = \kappa$. Then for every cardinal $\lambda$ with $\ink \leq \lambda \leq 2^\kappa$, there exists an almost disjoint family $\mathcal{A} \subseteq [\kappa]^\kappa$ of cardinality $\lambda$ such that $\Psi_\kappa(\mathcal{A})$ is not star-$\kappa$-Rothberger.
\end{Cor}

\bigskip
\noindent{\bf Acknowledgements:} The authors acknowledge the use of ChatGPT (OpenAI) for assistance with language editing and improving the clarity and presentation of the manuscript. The authors take full responsibility for the mathematical content and any remaining errors.

\medskip
\noindent{\bf Funding} No funding was received for conducting this study.

\medskip
\noindent{\bf Data Availability} No data was used for the research described in the article.

\medskip
\noindent{\bf Author Contributions} All authors contributed equally to this work.

\bigskip
\noindent{\bf Declarations}

\medskip
\noindent{\bf Conflict of Interest} All authors declare that they have no conflicts of interest.

\end{document}